%
%
%
%
%
%
\magnification=\magstep1   
\input amstex
\input pictex
\UseAMSsymbols
\hoffset=0truecm 
\hsize=125mm \vsize=185mm
\NoBlackBoxes
\parindent=8mm
\mathsurround=1pt
\font\gross=cmbx10 scaled\magstep1   \font\abs=cmcsc10
\font\rmk=cmr8  \font\itk=cmti8 \font\bfk=cmbx8   \font\ttk=cmtt8
\font\sfk=cmss7 
 \def\Ker{\operatorname{Ker}}
 \def\add{\operatorname{add}}
 \def\Aut{\operatorname{Aut}}
 \def\Im{\operatorname{Im}}
 \def\Ext{\operatorname{Ext}}
 
 \def\End{\operatorname{End}}
 \def\Rad{\operatorname{Rad}}
 \def\Soc{\operatorname{Soc}} 
 \def\Hom{\operatorname{Hom}}

 \def\can{\operatorname{can}}
 \def\Cok{\operatorname{Cok}}
 \def\Mono{\operatorname{Mono}}
 \def\Epi{\operatorname{Epi}}
 \def\Mimo{\operatorname{Mimo}}
 \def\Mepi{\operatorname{Mepi}}
 \def\incl{\operatorname{incl}}
 \def\mod{\operatorname{mod}}
 \def\I{\operatorname{I}}
 \def\P{\operatorname{P}}
 \def\R{{\text{\sfk R}}}
 \def\m{{\text{\bf m}}}
 \def\lto#1{\;\mathop{\longrightarrow}\limits^{#1}\;} 
 \def\lfrom#1{\;\mathop{\longleftarrow}\limits^{#1}\;} 
 \def\sto#1{\;\mathop{\to}\limits^{#1}\;}
 \def\lamod{\mod \Lambda}
 \def\undermod{\operatorname{\underline{mod}}}
 \def\overmod{\operatorname{\overline{mod}}}

 \def\underlamod{\undermod \Lambda}
 \def\overlamod{\overmod \Lambda}
 \def\hobject#1#2#3{{{\ssize #1\atop \ssize\phantom{#3}\downarrow #3}\atop
        \ssize #2}}

 \def\T#1{\qquad\text{#1}\qquad}
\def\E#1{{\parindent=1truecm\narrower\narrower\noindent{#1}\smallskip}}  
\def\qed{\phantom{m.} $\!\!\!\!\!\!\!\!$\nolinebreak\hfill\checkmark}
\headline{\ifnum\pageno=1\hfill %
    \else\ifodd\pageno \hfil\Rechts\hfil \else \hfil\Links\hfil \fi  \fi}
    \def\Links{\abs Ringel, Schmidmeier}
    \def\Rechts{\abs The Auslander-Reiten Translation in Submodule Categories}
\phantom{\noindent\rmk [rs-art-t, September 30, 2005]}
\bigskip\medskip

        \vglue1truecm
\centerline{\gross The Auslander-Reiten Translation}
        \smallskip
\centerline{\gross in Submodule Categories}
        \bigskip
\centerline{Claus Michael Ringel and Markus Schmidmeier}
        \bigskip
\centerline{\it Dedicated to Idun Reiten}
\centerline{\it on the occasion of her 60$^{\text{\it th}}$ birthday}
        \bigskip

\E{{\it Abstract.}
Let $\Lambda$ be an artin algebra or, more generally, a locally bounded
associative algebra, and $\Cal S(\Lambda)$ the category of all
embeddings $(A\subseteq B)$ where $B$ is a finitely generated 
$\Lambda$-module and $A$ is a submodule of $B$.
Then $\Cal S(\Lambda)$ is an
exact Krull-Schmidt category which has Auslander-Reiten sequences.
In this manuscript we show that the Auslander-Reiten translation
in $\Cal S(\Lambda)$ can be computed within $\lamod$ by using our
construction of minimal monomorphisms.  
If in addition $\Lambda$ is uniserial then any indecomposable nonprojective 
object in $\Cal S(\Lambda)$ is invariant under the sixth power of the 
Auslander-Reiten translation.  }

\medskip
\centerline{{\rmk 2000 Mathematics Subject Classification:}
{\itk Primary 16G70, Secondary 18E30}}
\centerline{{\rmk Keywords:} {\itk Auslander-Reiten sequences, approximations,
        triangulated categories}}


%
\bigskip
Let $\Lambda$ be an artin algebra, and $\lamod$ the category of finitely generated
$\Lambda$-modules (these are just the $\Lambda$-modules of finite length). 
The homomorphism category $\Cal H(\Lambda)$ has as objects the maps
$f$
in $\lamod$ and morphisms are given by commutative diagrams.
In this paper, we draw attention to the full subcategory $\Cal S(\Lambda)$ 
of $\Cal H(\Lambda)$ 
of all monomorphisms (or subobjects), but also to the corresponding subcategory 
$\Cal F(\Lambda)$ of $\Cal H(\Lambda)$ of all
epimorphisms (or factor objects).
Categories of the form $\Cal S(\Lambda)$ are much more complicated 
than the underlying module categories $\lamod$; for example,
if $\Lambda$ is a uniserial ring, then 
$\Lambda$ is of finite representation type, whereas the category 
$\Cal S(\Lambda)$ may have finitely or 
infinitely many indecomposable objects, or even 
be of wild representation type, depending on the Loewy length of $\Lambda$. 
Since Garrett Birkhoff in 1934 proposed the study of such 
submodule categories, 
they have proven to provide a rich source for classification problems,
and to attract the use of methods from various areas of algebra 
including representations of finite dimensional
algebras, lattices over tiled orders, representations of posets, and 
matrix classification.
In this manuscript we intend to lay the foundation for an Auslander-Reiten type 
theory of submodule categories.
Here is a preview:

\smallskip
If we consider a map $f\:A \to B$ of $\Lambda$-modules as an object of $\Cal H(\Lambda)$,
we will write either 
$$
 \big(A\sto fB\big) \quad\text{or}\quad \hobject ABf,
$$ 
but often also just $f$, 
whatever will be convenient and not misleading.
The category $\Cal H(\Lambda)$ is an abelian category, in fact it is 
equivalent to the category of finitely generated modules over the triangular matrix ring
$\big({\Lambda\atop0}{\Lambda\atop\Lambda}\big)$, and hence $\Cal S(\Lambda)$
as well as $\Cal F(\Lambda)$ are exact Krull-Schmidt category. 
We determine the projective and the injective objects
in $\Cal S(\Lambda)$ and $\Cal F(\Lambda)$.  For example, if $I$ is an
indecomposable injective $\Lambda$-module, then
$$\big(0\to I\big) \T{and} \big(I\sto{1}I\big)$$
are both indecomposable injective objects in $\Cal S(\Lambda)$, but clearly
the first is not injective in $\Cal H(\Lambda)$
(Proposition 1.4). 

\smallskip
In order to see that $\Cal S(\Lambda)$ has Auslander-Reiten sequences, we only have
to show that $\Cal S(\Lambda)$ is functorially finite in $\Cal H(\Lambda)$, according
to Auslander and Smal{\o} (Theorem 2.4 in [1]). But this is easy:
In the abelian category $\lamod$, every map $f\:A\to B$ can be factorized as 
the composition of an epimorphism $\Epi(f)$ 
and a monomorphism $\Mono(f)$. The factorization yields a morphism $f\to \Mono(f)$ in
$\Cal H(\Lambda)$ and this morphism is a left minimal $\Cal S(\Lambda)$-approximation for 
the object $f\in \Cal H(\Lambda)$.
To obtain a right minimal $\Cal S(\Lambda)$-approximation
for $f$, let $e'\:\Ker(f)\to \I\Ker(f)$ be an injective envelope and choose
an extension $e\:A\to \I\Ker(f)$ of $e'$.  We call the map
$$
 \Mimo(f) = [f\;e]\: \quad A\;\to\;B\oplus\I\Ker(f)
$$
a {\it minimal monomorphism} for $f$. In this way, we obtain a morphism 
$\Mimo(f)\to f$ in the category $\Cal H(\Lambda)$ and this morphism turns out to be 
the desired right minimal
$\Cal S(\Lambda)$-approximation for $f\in\Cal H(\Lambda)$ (Proposition 2.4).
Thus, the category $\Cal S(\Lambda)$ is functorially
finite in $\Cal H(\Lambda)$ and hence has Auslander-Reiten sequences. Note that 
$\Cal S(\Lambda)$-approximations provide recipes for calculating relative
Auslander-Reiten sequences in the subcategory $\Cal S(\Lambda)$
as soon as Auslander-Reiten sequences are known in $\Cal H(\Lambda)$.

\smallskip
In the module category $\Cal H(\Lambda)$, Auslander-Reiten translates are
computed via the Nakayama functor, as usual.  Using approximations
and the equivalence between $\Cal S(\Lambda)$ and $\Cal F(\Lambda)$ given
by the kernel and cokernel functors, corresponding Auslander-Reiten
sequences are obtained for the categories 
$\Cal S(\Lambda)$ and $\Cal F(\Lambda)$, see Proposition 3.2. (Surprisingly also the
converse is true:  Auslander-Reiten sequences in $\Cal S(\Lambda)$ 
and $\Cal F(\Lambda)$ give
rise to such sequences in the category $\Cal H(\Lambda)$, see
Proposition~3.5.)

\smallskip
Revisiting the construction of minimal monomorphisms, we show that
if two morphisms $f,g\:A\to B$ in $\lamod$ 
differ only by a map which factorizes through an injective $\Lambda$-module,
and if $B$ has no nonzero injective direct summands, then $\Mimo(f)$
and $\Mimo(g)$ are isomorphic objects in $\Cal S(\Lambda)$ (Proposition~4.1).

\smallskip
Our key result is this:
For an indecomposable nonprojective object 
$\big(A\sto fB\big)$ in 
$\Cal S(\Lambda)$, the Aus\-lan\-der-Reiten translate $\tau_{\Cal S}(f)$ 
can be computed directly within the category $\lamod$.  
This is done as follows: Let $g\:B\to C$ be the cokernel of $f$.  
Recall that the Auslander-Reiten
translation in $\lamod$ gives rise to a functor
$\tau_\Lambda\:\underlamod\to \overlamod$ where
$\underlamod$ and $\overlamod$ denote 
the factor categories of $\lamod$
modulo all maps which factorize through a projective or through an injective 
$\Lambda$-module, respectively. 
Next take a representative $h\:D\to E$ 
for the morphism $\tau_\Lambda(g)$
such that $D$ and $E$ have no nonzero injective direct summands.
As $h$ is determined uniquely,
up to a map which factorizes through an injective module, we will see that 
$\Mimo(h)$ is determined uniquely, up to 
isomorphism, as an object in $\Cal S(\Lambda)$.  
This is the Auslander-Reiten translate $\tau_{\Cal S}(f)$ for $f$
in the category $\Cal S(\Lambda)$ (Theorem~5.1).
Thus we may write: 
$$
 \tau_{\Cal S}(f)\;=\;\Mimo \tau_\Lambda \Cok(f)
$$

\smallskip
If $\Lambda$ is a self-injective algebra, the stable category
$\underlamod = \overlamod$ has the structure of a triangulated category. 
We observe that if $A\sto{\bar f}B\sto{\bar g}C\sto{\bar h}\Omega^{-1}A$
is a triangle (with $\Omega^{-1}$ the suspension functor)
and $f$ is a map representing $\bar f$, then
a map $g$ representing the rotate $\bar g$ is obtained as $g=\Cok\Mimo(f)$.
The functor $\tau_\Lambda$ 
commutes with $\Omega^{-1}$ and with
the rotation $\bar f\mapsto \bar g$ in a triangle,
and as a consequence, the 
formulae 
$$
 \overline{\tau_{\Cal S}^3(f)}\cong -\tau_\Lambda^3\Omega^{-1}(\bar f)
 \quad\text{and}\quad 
 \overline{\tau_{\Cal S}^6(f)}\cong \tau_\Lambda^6\Omega^{-2}(\bar f)
$$
hold (Theorem 6.2). 
In particular, if $\tau_\Lambda$ coincides with $\Omega^2$, it 
follows that 
$$
 \tau_{\Cal S}^3(f) \cong -\Mimo\Omega^{5}(f)
$$ 
for any indecomposable nonprojective object $f$ 
in $\Cal S(\Lambda)$ (Corollary 6.4).
For example, in the special case that $\Lambda$ is a commutative
uniserial ring, all the  functors $\tau_\Lambda,$ $\Omega^{2}$ 
and $\Omega^{-2}$ are 
equivalent to the identity functor on $\underlamod$ and then
the formula yields
$$
 \tau_{\Cal S}^6(f)\cong  f
$$ 
for $f$ an indecomposable nonprojective object
in $\Cal S(\Lambda)$ (Corollary 6.5).

\smallskip
Let $\big(C'\sto cC\big)$ be an indecomposable nonprojective object in
$\Cal S(\Lambda)$. The Aus\-lan\-der-Rei\-ten sequence 
$$\CD 0 @>>> \hobject{A'}A{a} @>f'>f> \hobject{B'}B{b} @>g'>g>
             \hobject{C'}C{c} @>>> 0 
\endCD $$
in $\Cal S(\Lambda)$ is made up from two short exact sequences in $\lamod$
given by the sequence $0\to A'\sto{f'}B'\sto{g'}C'\to 0$ of the 
submodules and the sequence $0\to A\sto fB\sto gC\to 0$ of the big modules.
By Proposition 7.2, these two sequences are ``usually'' split exact.
We list the exceptions and collect our findings about the structure
of the middle term $(b\:B'\to B)$.

	\smallskip
Let us stress that the categories
$\lamod$, $\Cal H(\Lambda)$ and $\Cal S(\Lambda)$ usually behave very differently. 
Consider for example the case of $\Lambda = \Lambda_n = 
k[T]/T^n.$ These $k$-algebras $\Lambda_n$ are all of
finite representation type: 
indeed, there is (up to isomorphism) precisely one indecomposable 
$\Lambda$-module of length $i$, for $1 \le i \le n$, and these are all the
indecomposables.
On the other hand, it is well-known that the matrix ring $\big({\Lambda_n\atop0}{\Lambda_n\atop\Lambda_n}\big)$ is representation finite only for $n \le 3$,
thus, for $n \ge 4$ there are infinitely many isomorphism classes
of indecomposable objects in $\Cal H(\Lambda_n)$, see for example [S,
Section~2]. For $n \le 5$,
the subcategory $\Cal S(\Lambda_n)$ of $\Cal H(\Lambda_n)$ consists 
of only finitely
many isomorphism classes of indecomposable objects, 
whereas there are infinitely many
isomorphism classes of indecomposable objects in $\Cal S(\Lambda_n)$, 
for any $n \ge 6,$ see [RS1].
One may consult these references and also [RS2] for proofs that the categories 
$\Cal H(\Lambda_4)$ and $\Cal S(\Lambda_6)$ are ``tame'', whereas the categories 
$\Cal H(\Lambda_n)$ 
for $n \ge 5$ and $\Cal S(\Lambda_n)$ for $n \ge 7$ are ``wild''.

\medskip\noindent{\it Notation:\/}
The condition that $\Lambda$ is an artin algebra can be weakened to 
requiring that $\Lambda$ is a locally bounded associative $k$-algebra or
a locally bounded $k$-spectroid [GR];
then also coverings of finite dimensional algebras are included.
We recall the corresponding definitions:

Let $k$ be a commutative local artinian ring and $\Lambda$ an 
associative $k$-algebra which need not have a unit
element, but it is required that $\Lambda$ equals the 
$k$-space $\Lambda^2$ of all
possible linear combinations of products in $\Lambda$. 
By $\lamod$ we denote the category of all $\Lambda$-modules $B$ which 
have finite length when considered as $k$-modules and which are 
{\it unitary\/} in the sense that $\Lambda B = B$ holds.
The algebra $\Lambda$ is said to be 
{\it locally bounded\/} if there is 
a complete set $\{e_i:i\in I\}$ of pairwise orthogonal primitive idempotents
such that each of the indecomposable projective 
modules $e_i\Lambda$ and $\Lambda e_i$ has finite length as a $k$-module, 
for $i\in I$. 
If $\Lambda$ is locally bounded such that 
each indecomposable projective $\Lambda$-module is injective
and 
each indecomposable injective $\Lambda$-module is projective, 
 then
we say that $\Lambda$ is a {\it self-injective algebra.}  An artin algebra $\Lambda$
is called {\it uniserial\/} if both modules $\Lambda_\Lambda$ and 
${}_\Lambda\Lambda$ have unique composition series. In this case, all
one-sided ideals are two-sided, namely the powers of the unique maximal
ideal $\m=\Rad\Lambda$.

For the terminology around almost split morphisms we refer the reader to [ARS] and [AS].
Here we use the term  ``Auslander-Reiten sequence'' for 
``almost split sequence''
and abbreviate ``left (right) minimal almost split map'' 
to ``source (sink) map''.
The Auslander-Reiten translation in a category $\Cal C$ is denoted by 
$\tau_{\Cal C}$, but in case $\Cal C = \mod \Lambda$, we write $\tau_\Lambda$. 

Finally, we want to apologize that our use of brackets 
when applying functions and
functors is not at all consistent. 
We have inserted brackets whenever we felt that this
improves the readability, 
but we have avoided multiple brackets whenever possible.

\medskip
This paper was written in 2001 as a general introduction to a proposed
volume devoted to the Birkhoff problem (dealing with subgroups of finite
abelian groups as well as with invariant subspaces of linear operators).
Unfortunately, we had to delay the Birkhoff project. Since the paper
seems of interest in its own right, we now have decided to publish it
independently. The authors are indebted to the referee,
and also to Aslak Bakke Buan and \O yvind Solberg, for helpful remarks
concerning the presentation of the results. 
	\bigskip

        \bigskip
\centerline{\bf  1\. Projective and Injective Objects}
        \medskip
In this section we determine the projective and the injective objects
in the categories $\Cal S(\Lambda)$ and $\Cal F(\Lambda)$ and their
associated sink and source maps. The injective objects in $\Cal S(\Lambda)$
are also called ``relatively injective'' or ``$\Ext$-injective'' as they
may not be injective in the category $\Cal H(\Lambda)$.
Dually, the projective objects in $\Cal F(\Lambda)$ may not be projective
when considered as objects in $\Cal H(\Lambda)$. 

\medskip
Let $\Lambda$ be an associative locally bounded $k$-algebra,
and let $U(\Lambda)=\big({\Lambda \atop 0}{\Lambda\atop\Lambda}\big)$
be the associative $k$-algebra of upper triangular matrices with 
coefficients in $\Lambda$.  First we recall well-known facts about 
$U(\Lambda)$ and about the category $\Cal H(\Lambda)$ of morphisms
between $\Lambda$-modules of finite $k$-length.

\medskip\noindent{\bf Lemma 1.1.} {\rm (Basic facts about $\Cal H(\Lambda)$)}
        {\it 
\item{1.} The $k$-algebra $U(\Lambda)$ is locally bounded.
\item{2.} The $k$-categories $\mod U(\Lambda)$ and $\Cal H(\Lambda)$ are
        equivalent.
\item{3.} $\Cal H(\Lambda)$ is an abelian Krull-Schmidt category.
\item{4.} Each object in $\Cal H(\Lambda)$ has a projective cover and an injective envelope.
\item{5.} The category $\Cal H(\Lambda)$ has Auslander-Reiten sequences.
\qed

}

\medskip The categories $\Cal S(\Lambda)$ and $\Cal F(\Lambda)$ are defined
to be the full subcategories of $\Cal H(\Lambda)$ which consist of all
objects $\big(A\sto fB\big)$ in $\Cal H(\Lambda)$ for which $f$ is a 
monomorphism or an epimorphism, respectively.  These three categories
are related by the kernel and cokernel functors.
$$\matrix 
\Cok: & \Cal H(\Lambda)\to\Cal F(\Lambda), & \quad (A\sto fB)
   \mapsto (B\lto{\text{can}}\Cok(f)),\cr
\Ker: & \Cal H(\Lambda)\to\Cal S(\Lambda), & \quad(B\sto gC)
   \mapsto(\Ker(g)\lto{\text{incl}}B).
\endmatrix$$

\medskip\noindent{\bf Lemma 1.2.} {\rm(Basic properties of $\Cal S(\Lambda)$ 
and $\Cal F(\Lambda)$)} {\it 
\item{1.} With the exact structure given by the category $\Cal H(\Lambda)$,
        the categories $\Cal S(\Lambda)$ and $\Cal F(\Lambda)$ are exact 
        Krull-Schmidt $k$-categories.
\item{2.} The category $\Cal S(\Lambda)$ is closed under kernels while
        $\Cal F(\Lambda)$  is closed under cokernels. Both categories
        are closed under extensions.
\item{3.} The restrictions of the kernel and cokernel functors
        $$\Ker\: \Cal F(\Lambda)\to \Cal S(\Lambda)\T{and}
                \Cok\:\Cal S(\Lambda)\to \Cal F(\Lambda)$$
        induce a pair of inverse equivalences.
\qed

}

\medskip The equivalence between $\Cal S(\Lambda)$ and $\Cal F(\Lambda)$
is useful to deduce the structure of the projective and the injective
objects in either of the two categories from the structure of the 
corresponding modules in $\Cal H(\Lambda)$ which are described in the following

\medskip\noindent
{\bf Lemma 1.3.} {\rm (Projective and injective modules in $\Cal H(\Lambda)$;
        the Nakayama functor)}
        {\it 
\parindent1cm
\smallskip
\item{P-1.} Let $P$ be an indecomposable projective $\Lambda$-module 
        with radical $\Rad P$. 
   The objects $(0\to P)$ and $(1_P\:P\to P)$ are indecomposable projective
   objects and have as sink maps the inclusions
   $$(0\to\Rad P)\to(0\to P)\quad\text{and}\quad
      (\Rad P\lto{\text{incl}}P)\to (P\sto1P),$$
   respectively.
\item{P-2.} Each indecomposable projective object arises in this way.

\smallskip
\item{I-1.} Let $I$ be an indecomposable injective $\Lambda$-module
         with socle $\Soc I$.
   The indecomposable injective objects $(I\to 0)$ and $(1_I\:I\to I)$ 
   have as source maps the canonical maps
   $$(I\to 0)\to (I/\Soc I\to 0)\quad\text{and}\quad
      (I\sto1I)\to(I\lto{\text{can}}I/\Soc I),$$
   respectively.
\item{I-2.} Each indecomposable injective object arises in this way.

\smallskip
\item{N-1.} The operation of the Nakayama functor $\nu_{\Cal H}$
           on $\Cal H(\Lambda)$ can be expressed in terms of the
           Nakayama functor $\nu=\nu_\Lambda$ in $\lamod$. For $P$ a
           projective $\Lambda$-module we have
           $$ \nu_{\Cal H}( 0\to P) \;=\; (\nu P\sto 1\nu P),\qquad
              \nu_{\Cal H}( P\sto 1P) \;=\; (\nu P\to 0).$$
\item{N-2.} On morphisms, the Nakayama functor $\nu_{\Cal H}$ is given in
           the obvious way, for example if $(0,f)\:(0\to P)\to (1_Q\:Q\to Q)$
           is a morphism between projective objects, then 
           $\nu_{\Cal H}(0,f)=
                (\nu f,0)\:(1_{\nu P}\:\nu P\to \nu P)\to (\nu Q\to 0)$.
\qed

}

\medskip
We can now describe the projective and the injective objects in
$\Cal S(\Lambda)$.

\medskip\noindent
{\bf Proposition 1.4.} 
        {\rm (Projective and injective objects in $\Cal S(\Lambda)$)} {\it 
\nopagebreak
\smallskip\parindent1cm
\item{P.} The projective objects in $\Cal S(\Lambda)$ and their 
        sink maps are as in Lemma 1.3-P.

\smallskip
\item{I-1.} Let $I$ be an indecomposable injective $\Lambda$-module.
        Then the map $(1_I\: I \to I)$ is an indecomposable 
   injective object in $\Cal S(\Lambda)$ and has  as
   source map the canonical map 
        $(1_I\: I\to I)\to(1_{I/\Soc I}\:I/\Soc I\to I/\Soc I)$.
\item{I-2.} If $I$ is an indecomposable injective $\Lambda$-module
        then $(0\to I)$ is 
   indecomposable (relatively) injective in $\Cal S(\Lambda)$ and 
   has as source map the inclusion 
   $(0\to I)\to (\incl\:\Soc I\to I)$.
\item{I-3.}  Each indecomposable (relatively) 
        injective object in $\Cal S(\Lambda)$ 
   arises in this way.

}

\smallskip\noindent
{\it Proof:}
The projective modules in $\Cal H(\Lambda)$ and their sink maps
are objects and morphisms in the category $\Cal S(\Lambda)$, hence
the statement in Lemma 1.3-P holds for $\Cal S(\Lambda)$. 
Similarly, the injective modules in $\Cal H(\Lambda)$ and their source
maps are in the category $\Cal F(\Lambda)$, so the statement in 
1.3-I holds for $\Cal F(\Lambda)$.  By applying the kernel functor
from 3.\ in Lemma 1.2 we obtain the objects and morphisms
in I-1 and I-2.  Since $\Cal H(\Lambda)$ has sufficiently many projective
and injective objects, so does $\Cal S(\Lambda)$. \qed

\medskip Let us add the dual statement for the category $\Cal F(\Lambda)$. 
Of course, we know
from 1.2 that the categories $\Cal S(\Lambda)$ and $\Cal F(\Lambda)$ are equivalent,
thus there is no intrinsic need to deal with both categories separately. On the other hand,
it may be useful for further references to have the precise formulations available.  

\medskip\noindent
{\bf Proposition 1.5.} 
        {\rm (Projective and injective objects in $\Cal F(\Lambda)$)} {\it 

\smallskip\parindent1cm
\item{P-1.} Let $P$ be an indecomposable projective $\Lambda$-module.
        Then the map $(1_I\: P \to P)$ is an indecomposable 
   projective object in $\Cal F(\Lambda)$ and has  as
   sink map the inclusion
        $(1_{\Rad P}\:\Rad P\to\Rad P)\to (1_P\:P\to P)$.
\item{P-2.} If $P$ is an indecomposable projective $\Lambda$-module
        then $(P\to 0)$ is 
   indecomposable (relatively) projective in $\Cal F(\Lambda)$ with 
        sink map the canonical map
        $(\can\:P\to P/\Rad P)\to (P\to 0)$.
\item{P-3.}  Each indecomposable (relatively) 
        projective object in $\Cal S(\Lambda)$ 
   arises in this way.

\smallskip
\item{I.} The injective objects in $\Cal F(\Lambda)$ and their 
        source maps are as in Lemma 1.3-I.
\qed

}

\medskip\noindent
{\it Example.\/} Consider the case that $\Lambda$ is uniserial
with maximal ideal $\m$ and Loewy length $n\geq 2$. 
In $\Cal S(\Lambda)$, there are two projective indecomposable objects
namely $P_1= (\Lambda = \Lambda)$ and $P_2= (0\subseteq \Lambda)$;
their sink maps are the inclusions
$$(\m\subseteq\Lambda)\quad\to\quad P_1\qquad\text{and}\qquad
         (0\subseteq \m)\quad\to\quad P_2,$$
respectively. 
One of them, the module $P_1=I_1$ is also injective, both in 
$\Cal H(\Lambda)$ and in $\Cal S(\Lambda)$, and has source map the
canonical map 
$$I_1\quad\to\quad (\Lambda/\m^{n-1}=\Lambda/\m^{n-1});$$
the second projective $P_2=I_2$ is relatively injective in
$\Cal S(\Lambda)$, its source map is given by the inclusion
$$I_2\quad\to\quad (\m^{n-1}\subseteq\Lambda).$$

       \bigskip
\centerline{\bf  2\. Left and Right Approximations}
\medskip
The categories $\Cal S(\Lambda)$ and $\Cal F(\Lambda)$ are functorially
finite in $\Cal H(\Lambda)$ and hence have Auslander-Reiten sequences.
To show this, we determine the left and the right approximation for each object
in $\Cal H(\Lambda)$, in each of the subcategories
$\Cal S(\Lambda)$ and $\Cal F(\Lambda)$.

\smallskip\noindent{\it Definitions:\/}
Let $\Cal S$ be a subcategory of a module category $\Cal C$,
and $C\in\Cal C$.  A morphism $f\:C\to S$ with $S\in\add\Cal S$ is a 
{\it left approximation\/} of $C$ in $\Cal S$ if the map 
$$\Hom(f,1_{S'})\:\quad \Hom_{\Cal C}(S,S')\to \Hom_{\Cal C}(C,S')$$
is onto for each $S'\in \Cal S$. Moreover, $f$ is {\it left minimal\/}
if each endomorphism $s\in\End_{\Cal C}(S)$ which satisfies
$fs=f$ is an isomorphism. A {\it minimal left approximation\/} is a 
left minimal left approximation.  Similarly, {\it minimal right
approximations\/} are defined.

\medskip
In the abelian category $\lamod$, 
a morphism $f\:A\to B$ has a factorization
$A\to\Im(f)\to B$ over the image. 
The two maps $f_1\:A\to \Im(f)$ and $f_2\:\Im(f)\to B$ 
are determined uniquely as objects in $\Cal H(\Lambda)$, up to isomorphism, 
and give rise to functors
$$ \matrix
\Epi: & \Cal H(\Lambda)\to\Cal F(\Lambda), 
   & \quad\big(A\sto fB\big)\mapsto\big(A\sto{f_1} \Im(f)\big),\cr
\Mono: & \Cal H(\Lambda)\to \Cal S(\Lambda), 
   & \quad\big(A\sto f B\big)\mapsto\big(\Im(f)\sto{f_2} B\big).
\endmatrix$$
There are functorial isomorphisms $\Epi\cong\Cok\Ker$ and $\Mono\cong\Ker\Cok$.

\medskip\noindent{\bf Lemma 2.1.} {\rm (Approximations given by
$\Mono$ and $\Epi$)} 

\noindent
{\it Let $(f\:A\to B)$ be an object in $\Cal H(\Lambda)$.

\smallskip
\item{1.} The map $(f_1,1_B)\:f\to \Mono(f)$ 
is a left minimal approximation of $f$ in $\Cal S(\Lambda)$.
\item{2.} The map $(1_A,f_2)\:\Epi(f)\to f$ 
is a right minimal approximation of $f$ in $\Cal F(\Lambda)$.

}

\smallskip\noindent The {\it proof\/} is immediate from the definitions. 
\qed

\medskip
For the object $(f\:A\to B)\in\Cal H(\Lambda)$
there is also a {\it right\/} minimal
approximation in $\Cal S(\Lambda)$ and a {\it left\/} minimal
approximation in $\Cal F(\Lambda)$, as we are going to show.

\medskip
First we define $\Mimo(f)$, the {\it minimal monomorphism} for $f$,
as follows: Let $e'\:\Ker(f)\to \I\Ker(f)$ 
be an injective envelope and choose an extension 
$e\:A \to \I\Ker(f)$
of $e'$, thus $e' = f'e$, where $f'\:\Ker(f)\to A$ is the
inclusion map. Then $\Mimo (f)$ is the map
$$
 \Mimo( f) = \bmatrix f & e \endbmatrix\: \quad 
                    A \;\to\; B\oplus \I\Ker(f).
$$

        \medskip\noindent
{\bf Lemma 2.2.} {\rm ($\Mimo$ is well-defined)}

\noindent{\it $\Mimo(f)$ is independent of the choice of $e$,
up to isomorphism in $\Cal H(\Lambda)$.}

\smallskip\noindent{\it Proof:}
Let $e_1, e_2\:A \to \I\Ker(f)$ be two extensions of $e'$, thus 
$e'=f'e_1=f'e_2$.
We see that the difference $e_2-e_1$ vanishes on $\Ker(f)$. Write 
$f=f_1f_2$ with $f_1$ an epimorphism and $f_2$ a monomorphism. 
Since $e_2-e_1$ vanishes on $\Ker(f)$, it factorizes 
through $f_1=\Cok (f')$, thus
$e_2-e_1 = f_1\hat e$ for some map $\hat e\:\Im(f)\to \I\Ker(f)$. 
Since $f_2$ is mono and the target 
$\I\Ker(f)$ of $\hat e$ is
injective, we can extend $\hat e$ to $B$: There is 
$\tilde e\:B \to \I\Ker(f)$ such that
$f_2\tilde e = \hat e$. 
Altogether we have $f \tilde e = f_1f_2\tilde e = f_1\hat e = e_2-e_1$. 
It follows that the representations in $\Cal S(\Lambda)$ given by 
$\bmatrix f & e_1\endbmatrix$ and by 
$\bmatrix f & e_2\endbmatrix$ are isomorphic.
$$
 \CD
 \ssize  A @>[f  \; e_1]>> \ssize B \oplus \I\Ker(f) \cr
 @|     @VV{\left[\smallmatrix 1 & \tilde e \cr 0 & 1\endsmallmatrix\right]} V\cr
 \ssize  A @>>[ f\; e_2]>  \ssize B\oplus \I\Ker(f) 
\endCD
$$ \qed

\medskip Dually one defines the 
{\it minimal epimorphism}, $\Mepi(f)$, for a map $f\:A\to B$ as follows.
The projective cover of the cokernel of $f$, $p'\:\P\Cok(f)\to \Cok(f)$
factorizes over the canonical map $f''\:B\to \Cok(f)$
so there is $p\:\P\Cok(f)\to B$ such that $p'=pf''$.
Let $\Mepi(f)$ denote the map
$$\Mepi(f)\;=\;\bmatrix f \cr p\endbmatrix\:\quad 
  A\oplus\P\Cok(f)\longrightarrow B;$$
then the dual version of the above result holds for $\Mepi$.

        \medskip\noindent
{\bf Lemma 2.3.} {\rm ($\Mepi$ is well-defined)}
\newline{\it $\Mepi(f)$ is independent of the choice of $p$,
up to isomorphism in $\Cal H(\Lambda)$.} \qed

\medskip
There are canonical maps $\Mimo(f)\to f$ and $f\to \Mepi(f)$,
$$\CD \hobject A{B\oplus I}{[f\;e]} @>1>\big[{1\atop0}\big]> 
                        \hobject ABf\qquad\text{and}
        \qquad\hobject ABf @>[1\;0]>1> 
                \hobject{A\oplus P}B{[f\;p]^t},
        \endCD$$
which give rise to a ``dual'' version of Lemma 2.1:

\medskip\noindent{\bf Proposition 2.4.} {\rm(Approximations defined by 
$\Mimo$ and $\Mepi$)} \newline {\it
Let $(f\:A\to B)$ be an object in $\Cal H(\Lambda)$.

\smallskip
\item{1.} The map $\Mimo(f)\to f$ 
is a right minimal approximation of $f$ in $\Cal S(\Lambda)$.
\item{2.} The map $f\to \Mepi(f)$ 
is a left minimal approximation of $f$ in $\Cal F(\Lambda)$.

}

\smallskip\noindent{\it Proof of the first statement:\/}
Let $(g\:C\to D)$ be an object in $\Cal S(\Lambda)$ and $(u,v)$ 
a morphism from $g$ to $f$.  We need to find $(u',v')\:g\to \Mimo(f)$
such that the composition with the above map 
$F\:\Mimo(f)\to f$ is just $(u,v)$.
Since $g\:C\to D$ is a monomorphism, the composition 
$C\sto uA\sto eI$ lifts to a map
$v_2\:D\to I$ such that $ue=gv_2$.  Then the pair $(u',v')=(u,[v,v_2])$ is a
morphism $g\to \Mimo(f)$ which satisfies the condition that
$(u',v')\,F=(u,v)$.
Thus, $F$ is a right approximation.  
It remains to show that $F$ is right minimal.
Let $(u,v)$ be an endomorphism of $\Mimo(f)$ such that $F=(u,v)\,F$ holds
and write the map $v\:B\oplus I\to B\oplus I$ as a matrix 
$v=\big[{v_{11}\atop v_{21}}{v_{12}\atop v_{22}}\big]$.
Since $F=(u,v)\,F$ holds, it follows that $u=1_A$, $v_{11}=1_B$, and 
$v_{21}=0$.  Now the condition that $(u,v)$ is a homomorphism amounts to
$[f\;e]\big[{1 \atop0} {v_{12}\atop v_{22}}\big]=
        1_A[f\;e]$, that is to say, $fv_{12}+ev_{22}=e$ must 
hold.  Restricting both sides to the kernel of $f$ (so that the composition
$fv_{12}$ vanishes) yields that the two maps $ev_{22},e\:\Ker(f)\to I$
are equal.  Since $I$ is the injective envelope of $\Ker(f)$, 
minimality implies that $v_{22}$ is an automorphism of $I$.  Hence 
$v=\big[{1_B\atop 0}{v_{12}\atop v_{22}}\big]$ is an automorphism 
of $B\oplus I$.  We have shown that $F$ is right minimal. \qed

\medskip As a consequence of Lemma 2.1 and Proposition 2.4 we obtain:

\medskip\noindent{\bf Theorem 2.5.} {\rm (Existence of Auslander-Reiten
sequences)} \newline{\it The subcategories $\Cal S(\Lambda)$
and $\Cal F(\Lambda)$ are functorially finite in $\Cal H(\Lambda)$ and
hence have Auslander-Reiten sequences.}

\smallskip\noindent{\it Proof:\/} According to the first statements
in Lemma 2.1 and Proposition 2.4, each object $f$ in $\Cal H(\Lambda)$
has a left and a right approximation in $\Cal S(\Lambda)$, this is to say,
$\Cal S(\Lambda)$ is functorially finite.  Similarly it follows from the 
second statements in Lemma 2.1 and Proposition 2.4 that
the category $\Cal F(\Lambda)$ is functorially finite, too.
According to [1, Theorem 2.4], those categories have Auslander-Reiten
sequences. \qed

       \bigskip
\centerline{\bf  3\. Transfer of Auslander-Reiten Sequences}
\nopagebreak\medskip
In this section we construct Auslander-Reiten sequences
for the categories $\Cal S(\Lambda)$ and $\Cal F(\Lambda)$ from
corresponding sequences in the module category $\Cal H(\Lambda)$.
Surprisingly, also the converse is possible:  Auslander-Reiten sequences in
$\Cal S(\Lambda)$ and $\Cal F(\Lambda)$ give rise to Auslander-Reiten
sequences in $\Cal H(\Lambda)$. 
First we need a lemma.

\medskip\noindent{\bf Lemma 3.1.} {\rm  (Kernels and Cokernels of 
Auslander-Reiten sequences)} \newline {\it
Suppose the following Auslander-Reiten sequence in 
the category $\Cal H(\Lambda)$ is given.
$$\CD 0 @>>> \hobject{A}{A_1}{a} @>f>f_1> 
        \hobject{B}{B_1}{b} @>g>g_1>
        \hobject{C}{C_1}{c} @>>> 0 
\endCD $$ 
If the kernel $\Ker(c)=(c'\:C'\to C)$ of the end term 
is not a projective object,
then the sequence obtained by applying the kernel
functor,
$$\CD 0 @>>> \hobject{A'}A{a'} @>f'>f> 
        \hobject{B'}B{b'} @>g'>g>
        \hobject{C'}C{c'} @>>> 0 
\endCD $$ 
is either split exact or almost split in  $\Cal S(\Lambda)$.  Dually, if 
the cokernel $\Cok(a)=(a''\:A_1\to A_1'')$ of the first term 
is not an injective object,
then the sequence obtained by using the cokernel functor,
$$\CD 0 @>>> \hobject{A_1}{A_1''}{a''} @>f_1>f_1''> 
        \hobject{B_1}{B_1''}{b''} @>g_1>g_1''>
        \hobject{C_1}{C_1''}{c''} @>>> 0 
\endCD $$ 
is either split exact or almost split
in $\Cal F(\Lambda)$.
}

\smallskip\noindent{\it Proof:\/} We only show the statement about
the sequence in $\Cal S(\Lambda)$, and for this sequence we only show
that the map $(g',g)\:(b'\:B'\to B)\to (c'\:C'\to C)$ is either a
split epimorphism or a right almost split morphism.  Note that this 
implies that $g'$ is onto.  

Let $(x'\:X'\to X)$ be an object in $\Cal S(\Lambda)$ and 
$(t',t)\:x'\to c'$ a morphism which is not a split epimorphism. 
Let $x\:X\to X_1$ be the cokernel for $x'$, factorize $c=c_1c_2$ as the
product of an epimorphism $c_1\:C\to \Im(c)$ and a monomorphism $c_2\:\Im(c)
\to C_1$, and let $t_1\:X_1\to \Im(c)$ be the cokernel map for $(t',t)$, 
this map satisfies
$xt_1=tc_1$.  Then $(t,t_1c_2)\:x\to c$ is not a split epimorphism, since
its kernel is not, and hence factorizes over the map $(g,g_1)$:
There is $(u,u_1)\:(x\:X\to X_1)\to (b\:B\to B_1)$ such that 
$(t,t_1c_2)=(u,u_1)(g,g_1)$. If $u'\:X'\to B'$ is the kernel map for $(u,u_1)$
then $(t',t)=(u',u)(g',g)$ factorizes over $(g',g)$. \qed

\medskip\noindent
{\bf Proposition 3.2.} {\rm (Transfer of AR sequences
        from $\Cal H(\Lambda)$ to 
        $\Cal S(\Lambda)$ and $\Cal F(\Lambda)$)} \newline {\it
Suppose that $0\to a\to b\to c\to 0$ is an Auslander-Reiten sequence
in $\Cal H(\Lambda)$.

\smallskip
\item{1.} If $c$ is an indecomposable nonprojective object in 
$\Cal F(\Lambda)$ then
$$0\to \Epi (a)\to \Epi (b)\to c\to 0$$
is an Auslander-Reiten sequence in $\Cal F(\Lambda)$ and
$$0\to \Ker (a)\to \Ker (b)\to \Ker (c)\to 0$$
is an Auslander-Reiten sequence in $\Cal S(\Lambda)$.

\smallskip\item{2.}
If $a$ is an indecomposable noninjective object in  $\Cal S(\Lambda)$ then
$$0\to a\to \Mono (b)\to \Mono (c)\to 0$$
is an Auslander-Reiten sequence in $\Cal S(\Lambda)$ and
$$0\to \Cok (a)\to \Cok (b)\to \Cok (c)\to 0$$
is an Auslander-Reiten sequence in $\Cal F(\Lambda)$.

}

\smallskip
\noindent{\it Proof:}\/
We show the first assertion.  Note that the
sequence $0\to \Epi(a)\to \Epi(b)\to c\to 0$ in $\Cal F(\Lambda)$
is the cokernel sequence
for $0\to \Ker(a)\to \Ker (b)\to \Ker (c)\to 0$ in $\Cal S(\Lambda)$ since 
$c\in\Cal F(\Lambda)$ satisfies $c=\Epi(c)=\Cok\Ker(c)$.
By Lemma~3.1, both sequences are either split exact or almost split
in their respective categories.  The sequences are not split exact since
the morphism $\Epi(b)\to c$ is the 
composition of the right approximation $\Epi(b)\to b$ and the right 
almost split  morphism $b\to c$, and hence is a right almost split morphism. 
\qed

\medskip\noindent
{\bf Corollary 3.3.} {\rm (The translation in $\Cal H(\Lambda)$,
$\Cal S(\Lambda)$ and $\Cal F(\Lambda)$)} {\it

\smallskip
\item{1.} If $c\in\Cal F(\Lambda)$ is indecomposable nonprojective then
   $\tau_{\Cal F}(c)=\Epi\tau_{\Cal H}(c).$
\item{2.} If $c\in\Cal S(\Lambda)$ is indecomposable nonprojective then
   $\tau_{\Cal S}(c)=\Ker\tau_{\Cal H}\Cok(c).$
\item{3.} If $a\in\Cal S(\Lambda)$ is indecomposable noninjective 
        then   $\tau^-_{\Cal S}(a)=\Mono\tau^-_{\Cal H}(a).$
\item{4.} If $a\in\Cal F(\Lambda)$ is indecomposable noninjective then
   $\tau^-_{\Cal F}(a)=\Cok\tau^-_{\Cal H}\Ker(a).$
\qed

}

\medskip In the remainder of this section we describe how Auslander-Reiten
sequences in $\Cal S(\Lambda)$ and $\Cal F(\Lambda)$ give rise to 
Auslander-Reiten sequences in $\Cal H(\Lambda)$. 

\medskip\noindent
{\bf Lemma 3.4.} {\rm (The functors $\Epi$ and $\Mono$ reflect some
split morphisms)} \newline {\it
Let $(f,g)\:\big(A\sto aB\big)\to \big(C\sto cD\big)$ 
be a morphism in $\Cal H(\Lambda)$
and let $h\: \Im (a)\to \Im (c)$ be the induced map on the images.
$$\CD
  \ssize A  @>a_1>>  \ssize \Im(a)  @>a_2>>  \ssize B\cr
 @VfVV  @VhVV  @VVgV \cr
  \ssize C @>>c_1>  \ssize \Im(c)  @>>c_2> \ssize D
\endCD$$
\item{1.} If $(f,g)$ is split monic  (split epic)
then $\Epi(f,g)=(f,h)$ and $\Mono(f,g)=(h,g)$
are both split monic (split epic).

\item{2.} If $a_2\:\Im(a)\to B$ is an injective envelope,
and if $\Epi(f,g)$ is split monic, then so is $(f,g)$.

\item{3.} If $c_1\:C\to\Im(c)$ is a projective cover, 
and if $\Mono(f,g)$ is split epic, then so is $(f,g)$. 

}

\smallskip\noindent
{\it Proof:}\/ The first assertion is clear since $\Epi$ and $\Mono$ are
functors; the third assertion is dual to the second, so we only consider
the second assertion.  Suppose that $(f,h)$ is split monic,
so there are maps $u\:C\to A$ and $w\: \Im(c)\to \Im(a)$ such that
$ua_1=c_1w$, $fu=1_A$, $hw=1_{\Im (a)}$. Since $B$
is injective, the map $wa_2$ extends to $D$:
There is $v\:D\to B$ such that $c_2v=wa_2$.
Since $a_2gv=hc_2v=hw a_2=a_2$,
the left minimality of the injective envelope $a_2$ implies that $gv$ is an 
automorphism of $B$. It follows that the morphism $(f,g)$ is split monic.
\qed

\medskip\noindent
{\bf Proposition 3.5.} {\rm (Transfer of AR sequences 
        from $\Cal S(\Lambda)$ and 
        $\Cal F(\Lambda)$ to $\Cal H(\Lambda)$)} {\it
\nopagebreak\smallskip
\item{1.}
Let $(c\:C\to C'')$ be an indecomposable nonprojective 
object in $\Cal F(\Lambda)$.
If the first two rows in the diagram
$$\CD
 \ssize 0  @>>> \ssize A @>f>> \ssize B @>>>\ssize C @>>> \ssize 0 \cr
     @.          @VaVV        @VVbV       @VVcV \cr
 \ssize 0 @>>>  \ssize A'' @>f''>> \ssize B'' @>>> \ssize C'' @>>> \ssize 0 \cr
     @.          @V\iota VV          @VVdV  @|       \cr
 \ssize 0  @>>> \ssize \I(A'') @>g>> \ssize D @>>>\ssize C'' 
                                                    @>>> \ssize 0\cr
\endCD$$
define an Auslander-Reiten sequence in $\Cal F(\Lambda)$ and if the third row 
is obtained as push-out along the injective envelope $\iota:A''\to \I(A'')$
then the first row and the third row define an Auslander-Reiten sequence
in $\Cal H(\Lambda)$.
\item{2.}Let $(a\:A'\to A)$ be an indecomposable noninjective
object in $\Cal S(\Lambda)$.
If the two lower rows in the diagram
$$\CD
 \ssize 0 @>>> \ssize A' @>>> \ssize D @>>h> \ssize \P(C') @>>> 
                                                            \ssize 0 \cr
     @.          @|         @VVdV  @VV\pi V       \cr
 \ssize 0 @>>>  \ssize A' @>>> \ssize B' @>>g'> \ssize C' @>>> \ssize 0 \cr
     @.          @VaVV        @VVbV       @VVcV \cr
 \ssize 0  @>>> \ssize A @>>f> \ssize B @>>>\ssize C @>>> \ssize 0 \cr
\endCD$$
define an Auslander-Reiten sequence in $\Cal S(\Lambda)$ and if the first row 
is obtained as pull-back along the projective cover  $\pi\:\P(C')\to C'$
then the first row and the third row define an Auslander-Reiten sequence
in $\Cal H(\Lambda)$.

}

\smallskip\noindent{\it Proof of the first statement:}\/
The map $(f,g)$ in $\Cal H(\Lambda)$ is not a split monomorphism
since $\Epi(f,g)=(f,f'')$ is not, by Lemma 3.1.
We show that $(f,g)$ is left almost split. Suppose
$$(r,s): \quad \big(A\lto{a\iota}\I(A'')\big)\quad\longrightarrow
   \quad \big(X\sto xY\big)$$
is a morphism in $\Cal H(\Lambda)$ which is not a split monomorphism. 
Factorize $x$ over the image as $X\sto{x_1} \Im(x)\sto{x_2}Y$, so
$x=x_1x_2$.
Then there is $t:A''\to \Im(x)$ such that the two squares commute.
$$\CD \ssize A @>r>> \ssize X \cr
@VaVV @VVx_1V \cr
\ssize A'' @>t>> \ssize \Im(x)\cr
@V\iota VV  @VVx_2 V \cr
\ssize \I(A'') @>s>> \ssize Y
\endCD $$
By Lemma 3.4, part 2, the morphism 
$(r,t)$ is not split monic in $\Cal F(\Lambda)$ and hence factorizes over 
$(f,f'')$: There are $v:B\to X$, $v'':B''\to \Im(x)$ such that 
$$bv''=vx_1, \quad\text{and}\quad (r,t)=(f,f'')(v,v'').$$
Then $f''(v''x_2)= tx_2 = \iota s$, so we obtain 
$w:D\to Y$ such that $dw=v''x_2$
and $gw=s$, by the push-out property for $Y$.
Thus, $(v,w)$ is a morphism in $\Cal H(\Lambda)$ and 
our test map $(r,s)$ factorizes over $(f,g)$:
$(r,s)=(f,g)(v,w)$.
For the proof that $(f,g)$ is a source map in 
the module category $\Cal H(\Lambda)$, it remains
to check that $(f,g)$ is left minimal, and this follows from the
indecomposability of the cokernel $(C\to C'')$. 
\qed

       \bigskip
\centerline{\bf  4. Minimal Monomorphisms and the Stable Category}
        \medskip
Returning to the investigation of minimal monomorphisms, we show that
if $f,g\:A\to B$ are two maps which differ by a morphism which factorizes 
through an injective module, then $\Mimo(f)$ and $\Mimo(g)$ are isomorphic
as objects in $\Cal S(\Lambda)$.  

\medskip
First we verify three claims.

        \medskip\noindent
{\bf  Claim 1.} {\it 
Let $f,g\:A \to B$ be maps in $\lamod$ such that $g-f$ factorizes through an
injective module. Let $h\:A \to \I(A)$ be an injective envelope.
Then the objects in $\Cal S(\Lambda)$ 
given by the monomorphisms $[f\; h]$ and $[g\; h]$ are isomorphic.} 

        \smallskip\noindent
{\it Proof:} The map $g-f$ factorizes through the injective envelope $h$, 
so there is
$u\:\I(A)\to B$ such that $g-f=hu$.
The following commutative diagram shows that $[f\;h]$ and $[g\;h]$ 
are isomorphic objects in $\Cal H(\Lambda)$. 
$$
\CD
 \ssize A @>[f  \; h]>> \ssize B \oplus \I(A) \cr
 @|     @VV{\left[\smallmatrix 1 & 0 \cr u & 1\endsmallmatrix\right]} V\cr
 \ssize A @>[g  \; h]>> \ssize B \oplus \I(A)
\endCD
$$
\qed

        \medskip\noindent
{\bf Claim 2.} {\it 
Given $f:A\to B$ and a map $h:A\to I$ with $I$ injective 
such that $[f\; h]\: A \to B \oplus I$ is a monomorphism, 
then there is an injective module $J$ and a commutative diagram 
$$
\CD 
 \ssize A @>[ f \; e \; 0]>> \ssize B\oplus \I\Ker (f) \oplus J \cr
 @|     @VV d V\cr
 \ssize A @>[ f \; h]>> \ssize B\oplus I
\endCD
$$
with $d$ an isomorphism and $e\:A\to \I\Ker(f)$ an extension of
the injective envelope $\Ker(f) \to \I\Ker(f).$}

\smallskip
The diagram shows that the maps given by the two rows are isomorphic as objects in
$\Cal H(\Lambda)$. 
Note that in $\Cal H(\Lambda)$, the object given by the upper row decomposes as the
direct sum of the two objects 
$([f\; e]\:A \to B\oplus \I\Ker(f))$ and 
$(0 \to J)$. Up to isomorphism, the first one is just $\Mimo(f)$. 
        
\smallskip\noindent
{\it Proof:\/}  
Since $[f \; h]$ is a monomorphism, the restriction of $h$ to
$\Ker(f)$ is injective, thus $I$ contains a submodule isomorphic to $\Ker(f)$
and therefore $I$ decomposes as 
$I = \I\Ker(f) \oplus J$ for some injective $\Lambda$-module $J$.
Then $h$ has the form $h = [ e \; h_2 ]$, where 
$ e\: A \to \I\Ker(f)$ is an extension of the inclusion 
$\Ker(f)\to \I\Ker(f)$,
and $ h_2\:A \to J$ satisfies $h_2\Ker(f) = 0.$ 
Write $f = f_1f_2$ with $f_1$ an epimorphism and $f_2$ a monomorphism.
Since $h_2$ vanishes on $\Ker(f)$, we can factorize 
$h_2$ through $f_1$, say $h_2 = f_1v$ for some map 
$v\:\Im(f)\to J$. Since $J$ is injective and $f_2$ 
a monomorphism, we obtain a lifting 
$w\:B\to J$ of 
$v$ to $B$, thus $v = f_2w$. 
Altogether we see that 
$h_2 = f_1v = f_1 f_2 w = fw$. 
This shows that the following diagram commutes.
$$
\CD 
 \ssize A @>[ f \; e \; 0]>> \ssize B\oplus \I\Ker (f) \oplus J \cr
 @|     @VV{\left[\smallmatrix 1 & 0 &w\cr 0 & 1 & 0\cr 0 & 0 & 1 \endsmallmatrix\right]} V\cr
 \ssize A @>[ f \; e \; h_2]>> \ssize B\oplus \I\Ker(f) \oplus J
\endCD
$$
\qed

\medskip
\noindent
{\bf Claim 3.} {\it
Suppose that $B$ has no nonzero injective direct summands
and that $f,g\:A\to B$ are maps such that 
$g-f$ factorizes through an injective
object. Then $\I\Ker(f)=\I\Ker(g)$.
}

\smallskip\noindent
{\it Proof:\/} (a) First we show that if a map $h:A \to B$
factorizes through an injective object, say $h=h_1h_2$ where 
$h_1\:A\to I$ and $h_2\:I\to B $, then $\Soc\Ker(h)=\Soc (A)$.
Indeed, if there were a simple submodule $S$ of $A$ such that the composition
$$S\lto{\text{incl}} A \lto{h_1} I \lto{h_2} B $$
is nonzero then we 
would obtain that $B$ has a nonzero injective direct summand ---
in contradiction to our assumption on $B$. 

\smallskip\noindent
(b) As a consequence we obtain that
if $g-f$ factorizes through an injective object,
then $\Soc\Ker(f)=\Soc\Ker(g)$ holds. By (a), $\Soc(A)\subseteq \Ker (g-f)$
and hence
$$\Ker(f)\cap \Soc(A)\; = \; \Ker(f+(g-f))\cap \Soc(A) \;=\; 
  \Ker(g) \cap \Soc (A).$$
Thus, $\I\Ker(f)=\I\Soc\Ker(f) = \I\Soc\Ker(g) = \I\Ker(g)$ holds.
\qed

\medskip\noindent
{\bf Proposition 4.1.} {\rm ($\Mimo$ independent of maps which 
        factorize through injective)}

\noindent{\it 
Suppose that $B$ has no nonzero injective direct summands.
Let $f,g\:A\to B$ be maps in $\lamod$ such that
$g-f$ factorizes through an injective $\Lambda$-module. 
Then the objects $\Mimo(f)$, $\Mimo (g)$ are isomorphic in $\Cal S(\Lambda)$. 
}

\smallskip\noindent
{\it Proof:\/} Let $h\:A\to \I(A)$ be an injective envelope,
then the objects $[f \; h]$, $[g \; h]$ are isomorphic by Claim 1.
According to Claim 2, there exist injective modules $J_1$, $J_2$,
extensions $e_1\:A\to \I\Ker(f)$, $e_2\:A\to \I\Ker(g)$ of the
inclusion maps $\Ker(f)\to \I\Ker(f)$, $\Ker(g)\to \I\Ker(g)$, 
respectively, and isomorphisms $d_1$, $d_2$ such that the diagram below
is commutative.
$$\CD
 \ssize A @>>[f \; e_1 \; 0]> \ssize B\oplus\I\Ker(f)\oplus J_1 \cr
 @|          @VVd_1V \cr
 \ssize A @>>[f \; h ]> \ssize B\oplus \I(A) \cr
 @|          @VV\cong V \cr
 \ssize A @>>[g \; h ]> \ssize B  \oplus \I(A) \cr
 @|          @AAd_2A \cr
 \ssize A @>>[g \; e_2 \; 0]> \ssize B \oplus\I\Ker(g)\oplus J_2
\endCD$$
According to Claim 3, the $\Lambda$-modules $\I\Ker(f)$, $\I\Ker(g)$ 
are isomorphic. By the Krull-Remak-Schmidt theorem for $\Lambda$-modules, 
$J_1\cong J_2$
follows. Note that the top row in the diagram when considered as an
object in $\Cal H(\Lambda)$ is isomorphic to the direct sum
$\Mimo(f)\oplus (0\to J_1)$ while the bottom row is isomorphic to
$\Mimo(g)\oplus (0\to J_2)$. Applying the Krull-Remak-Schmidt theorem in the 
category $\Cal H(\Lambda)$ we obtain that $\Mimo(f)\cong \Mimo (g)$.
\qed

\medskip
The following example shows that the condition that $B$ has no 
nonzero injective direct summands cannot be omitted.

\smallskip\noindent
{\it Example.} Let $\Lambda$ be a uniserial ring of 
Loewy length 2 with $m$ a
generator of the maximal ideal $\m$.
Denote by $\mu_m\:\Lambda\to\Lambda$ the multiplication by $m$.
With $f=1\:\Lambda\to\Lambda$ and $g=\mu_m\:\Lambda\to \Lambda$,
clearly $g-f$ factorizes through an injective $\Lambda$-module
but $\Mimo(f)=f$ and $\Mimo(g)=\big(\Lambda\lto{[\mu_2\;1]}\Lambda
\oplus \Lambda\big)$ are not isomorphic as $\Mimo(g)\cong f\oplus 
(0\to\Lambda)$.

\medskip
The criterion in Proposition 4.1 for $\Mimo(f)\cong\Mimo(g)$ can be refined
to obtain an equivalent condition.

\medskip\noindent
{\bf Theorem 4.2.} {\rm (An equivalent condition for $\Mimo(f)\cong\Mimo(g)$)}
\newline {\it 
Let $f,g\:A\to B$  be maps in $\lamod$
and suppose that $B$ has no nonzero injective direct summands.
Then $\Mimo(f)$ and $\Mimo(g)$ are isomorphic in $\Cal S(\Lambda)$ 
if and only if there are 
automorphisms $a$ of $A$ and $b$ of $B$ such that
$fb-ag$ factorizes through an injective $\Lambda$-module.}

\medskip\noindent
{\it Proof:} Assume first that $\Mimo(f)\:A\lto{[f\;e_1]}B\oplus I$
and $\Mimo(g)\:A\lto{[g\;e_2]}B\oplus I$ are isomorphic.
Then there are maps $a\in\Aut A$ and 
$h=\Big({b \atop h_{12}}\;{h_{21}\atop h_{22}}\Big)
\in\Aut (B\oplus I)$
such that $[f\;e_1]h=a[g\;e_2]$.  Thus, 
$ag=fb+e_1h_{12}$, so $ag-fb$ 
factorizes through the injective $\Lambda$-module $I$.
Moreover, as $B$ and $I$ have no indecomposable direct summands
in common, it follows that  the map $h$ is an automorphism if and only if
both $b$ and $h_{22}$ are automorphisms.

For the converse, assume that $fb-ag$ factorizes through
an injective $\Lambda$-module.
By Proposition 4.1, the objects $\Mimo(fb)$ and $\Mimo(ag)$
are isomorphic; assume they are given by maps 
$$\Mimo(fb)=[fb\;\;e_1], \;
\Mimo(ag)=[ag\;\;e_2]\: 
\qquad A\longrightarrow B\oplus I.$$
Thus, $\Mimo(f)\cong\Mimo(g)$, as indicated by
the following diagram.
$$ \CD
\Mimo(f)\:\qquad @. \ssize A @>[f\;e_1]>> 
                        \ssize B\oplus I \cr
@.  @|  @VV\Big[{b\atop 0}\;{0\atop 1}\Big]V \cr
\Mimo(fb)\:\qquad @. \ssize A @>[fb\;\;e_1]>> 
                        \ssize B\oplus I \cr
@.  @|  @VV\cong V \cr
\Mimo(ag)\:\qquad @. \ssize A @>[ag\;\;e_2]>>
                        \ssize B\oplus I \cr
@.  @Va VV @| \cr
\Mimo(g)\:\qquad @. \ssize A @>[g\;\;a^{-1}e_2]>>
                        \ssize B\oplus I \cr
\endCD
$$
For the last step note that if $e_2\:A\to I$
is an extension to $A$ of an injective envelope for $\Ker(ag)$,
then since $\Ker(ag)=a^{-1}\Ker(g)$, one obtains that 
$a^{-1}e_2$ is an extension to $A$ of an 
injective envelope for $\Ker (g)$.  By Lemma 2.2, the object
$([g\;\;a^{-1}e_2]\:A\to B\oplus I)$ is $\Mimo(g)$, up to isomorphism. 
\qed

\medskip
There is the following dual result for minimal epimorphisms.

\medskip
\noindent{\bf Theorem 4.3.} {\rm (An equivalent condition for 
        $\Mepi(f)\cong\Mepi(g)$)}
\newline {\it Let $f,g\:A\to B$  be maps in $\lamod$
and suppose that $A$ has no nonzero projective direct summands.
Then $\Mepi(f)$ and $\Mepi(g)$ are isomorphic 
in $\Cal F(\Lambda)$ if and only if there are 
automorphisms $a$ of $A$ and $b$ of $B$ such that
$fb-ag$ factorizes through a projective $\Lambda$-module.}\qed

       \bigskip
\centerline{\bf  5. The Auslander-Reiten Translation}
        \medskip

In Chapter 3 we have seen how the Auslander-Reiten translations
in the categories $\Cal H(\Lambda)$, $\Cal S(\Lambda)$, and $\Cal F(\Lambda)$
are related. 
Using this, we develop a formula to compute the Auslander-Reiten translate
$\tau_{\Cal S}(f)$  of an object $(f\:A\to B)$ in $\Cal S(\Lambda)$ 
directly in the category $\lamod$.  Here is the main result:

        \medskip\noindent
{\bf Theorem 5.1.} {\rm (The Auslander-Reiten translation
        in $\Cal S(\Lambda)$)} \newline
{\it
For an object $\big(A\sto fB\big)$ in $\Cal S(\Lambda)$, the Auslander-Reiten
translate is given by}
$$
 \tau_{\Cal S}(f) = \Mimo\tau_\Lambda\Cok(f).
$$
        \medskip
Note that this means the following: 
We start with the monomorphism $f$ and form 
its cokernel $g$. 
We apply $\tau_\Lambda$ to this map. 
Recall that $\tau_\Lambda(g)$ is only well-defined in the category 
$\overlamod$ (obtained from $\lamod$ by factoring out 
all the maps which factorize through an injective object). 
Represent $\tau_\Lambda(g)$ by a morphism $h\:D\to E$ in $\lamod$ 
such that $D$ and $E$ have no nonzero injective direct summands.
Now apply $\Mimo.$ As we have seen in Chapter 4,
this yields a well-defined object in $\Cal S(\Lambda)$, up to isomorphism.

        \medskip\noindent
{\it Proof:\/}
We proceed as follows.  
Let $(f\:A\to B)$ be an object in $\Cal S(\Lambda)$ with 
cokernel $(g\:B\to C)$. 
In Step 1 we obtain an exact sequence
defining $\tau_\Lambda(g)$; in Step 2 we construct 
an exact sequence
involving $\tau_{\Cal H}(g)$, from which $\tau_{\Cal S}(f)$ is computed
as $\Ker\tau_{\Cal H}(g)$ by Corollary 3.3.
The formula for the computation of a kernel in Step 3 is used
in Step 4 to relate the two exact sequences, and finally, in Step 5, we 
obtain our result.

\smallskip\noindent{\it Step 1}.
For an epimorphism $g\: B\to C$ in $\lamod$ we determine $\tau_\Lambda(g)$.
We start with a minimal projective presentation of $B$, say
$$ Q\;\lto d \; P \;\lto e \;B\;\longrightarrow\; 0$$
and a projective presentation of $C$ using the map $eg\:P\to C$
and a projective cover $t\: R\to \Ker(eg)$ to obtain the
following commutative diagram with exact rows.
$$\CD
 \ssize Q     @>d>>   \ssize P     @>e>>  \ssize B   @>>>  \ssize 0    \cr
 @VsVV                @VV1V               @VVgV                      \cr
 \ssize R     @>>t>   \ssize P     @>>eg>  \ssize C @>>> \ssize 0  \cr
\endCD $$
By applying the Nakayama functor 
$\nu_\Lambda=\nu=D\Hom_\Lambda(-,{_\Lambda\Lambda})$ in $\lamod$ to the left
square we arrive at the diagram defining $\tau_\Lambda(g)$:
$$\CD
 \ssize 0   @>>> \ssize  \tau_\Lambda B @>v>>   \ssize \nu Q   @>\nu d >>
                                                \ssize \nu P           \cr
 @.              @V\tau_\Lambda(g)VV          @VV \nu s V  @VV 1V     \cr
 \ssize 0   @>>> \ssize \tau_\Lambda C @>>w> \ssize \nu R  @>>\nu t >
                                                \ssize \nu P      
\endCD $$
(Since $eg$ is not necessarily a projective cover, $\nu t$ is not
necessarily an injective envelope.)

\smallskip\noindent{\it Step 2}.
From a projective presentation for $g$ in the category $\Cal H(\Lambda)$ 
we construct an exact sequence which involves
$\tau_{\Cal H}(g)$:  Being an epimorphism,
the object $g$ has a projective cover given by
$(1_P\: P\to P)$ with $P$ as in Step 1. Using the maps $d$ and $t$
we obtain the following projective presentation for $g$ in $\Cal H(\Lambda)$.
$$\CD
\hobject QQ1 \oplus \hobject 0R{ } 
      @>{\left[{d\atop 0}\right]}>{\left[{d\atop t}\right]}> \hobject PP1 
      @>e>eg> \hobject B{C}{g} @>>> 0 
\endCD $$
Note that this differs in general from a minimal projective presentation;
in order to obtain a minimal one, we would have to split off a direct 
summand of the form $(0\to S)$ where $S$ is a direct summand of $R$
(and $S$ is also a direct summand of $Q$).
Applying the Nakayama functor $\nu_{\Cal H}$ (see Lemma 1.3-N) to the 
morphism between the projective modules, we obtain
the following sequence
$$\CD
0 @>>> \tau_{\Cal H}\Big( \!\!\!\hobject B{C}{g}\!\!\Big)
\oplus \hobject{\nu S}{\nu S}1 @>>> \hobject{\nu Q}0{} \oplus
\hobject{\nu R}{\nu R}1 @>{\left[{\nu d\atop \nu t}\right]}>0> 
  \hobject{\nu P}0{}
\endCD$$
in which the additional projective summand $(0\to S)$ gives rise to the 
injective direct summand $(1_{\nu S}\: \nu S\to \nu S)$.

\smallskip\noindent{\it Step 3}.
For the computation of the kernel of a map $u\: U\to V$ 
we observe that if 
$$\CD 0 @>>> \hobject{U'}0{} @>u'>0> \hobject U{V}u @>>>
                    \hobject {W'}Ww \endCD $$
is an exact sequence in $\Cal H(\Lambda)$ with $(w\:W'\to W)$
in $\Cal S(\Lambda)$, then $\Ker u=u'$. 

\smallskip\noindent{\it Step 4}.
In the following diagram, the right two columns are split exact and
form a commutative diagram in $\Cal H(\Lambda)$; the left column is just
the kernel sequence.  Note that the sequence in Step 2 
involving $\tau_{\Cal H}(g)$ 
occurs as the middle row, while the sequence in Step~1 defining 
$\tau_\Lambda B$ occurs as the sequence of source modules in the top row. 
$$\CD
 @.  0  @.  0  @.  0 \cr
 @.  @VVV          @VVV          @VVV     \cr
 0 @>>> \hobject{\tau_\Lambda B}0{}  @>v>0>  \hobject{\nu Q}0{}
         @>\nu d>0>  \hobject{\nu P}0{}   \cr
 @.  @VV{\gamma \atop 0}V  @VV{ [1 \; 0]\atop 0}V  @VV{ 1\atop 0 }V \cr
 0 @>>>  \tau_{\Cal H}\Big( \!\!\!\hobject B{C}{g}\!\!\Big)
       \oplus \hobject{\nu S}{\nu S}1 
   @>>> \hobject{\nu Q}0{}\oplus \hobject{\nu R}{\nu R}1 
   @>{\left[{\nu d\atop \nu t}\right]}>0>
   \hobject{\nu P}0{} \cr
 @. @VVV  @VV{\left[{0\atop 1}\right]\atop\left[{0\atop 1}\right]}V @VVV \cr
 0  @>>> \hobject {\nu R}{\nu R}1 @>1>1> \hobject {\nu R}{\nu R}1 @>>> 0  \cr
 @.  @.  @VVV @VVV \cr
 @. @. 0 @. 0 \cr
\endCD $$
Since the object $1_{\nu R}\:\nu R\to \nu R$ is in $\Cal S(\Lambda)$,
we can use the left hand column for the computation of
$\Ker \tau_{\Cal H}(g)$, using Step 3:
$$\Ker\left(\tau_{\Cal H}(g)\oplus (1_{\nu S}\:\nu S\to \nu S)\right)=\gamma.$$
In order to identify $\gamma$, only the first entries in the above
diagram play a role, thus we have to deal with the following diagram.
$$(*)\qquad \CD
  \ssize 0  @>>> \ssize \tau_\Lambda B @>v>>  \ssize \nu Q @>\nu d>> 
                                                        \ssize \nu P  \cr
  @.  @VV\gamma V  @VV[1\;0]V  @| \cr
  \ssize 0  @>>> \ssize X  @>>>  \ssize \nu Q\oplus \nu R 
                                @>>\left[{\nu d\atop \nu t}\right]>\ssize\nu P
\endCD$$
where $X$ is the first component of 
$\tau_{\Cal H}(g)\oplus(1_{\nu S}\:\nu S\to \nu S)$.
Compare this with the following diagram.
$$(**)\qquad\CD
 \ssize 0 @>>>\ssize \tau_\Lambda B @>v>> \ssize \nu Q @>\nu d>>\ssize\nu P \cr
 @.              @VV[\tau_\Lambda(g)\;v]V   @VV[1\;0]V            @|     \cr
 \ssize 0 @>>> \ssize \tau_\Lambda C\oplus \nu Q
          @>>{\left[{0\atop 1}{w\atop -\nu s}\right]}>
    \ssize\nu Q\oplus \nu R @>>{\left[{\nu d\atop \nu t}\right]}> \ssize \nu P
\endCD $$
A glance back to Step 1 shows that the diagram is commutative and has an
exact upper row.  It only has to be confirmed that the lower row is exact.
Clearly, the $2\times 2$-matrix is a monomorphism, and the composition
of the last two maps is zero. 
Now, given $(x,y)\in\nu Q\oplus \nu R$ with $(x)\nu d+(y)\nu t = 0$,
then $\nu d=\nu s\,\nu t$ yields that $((x)\nu s+ y)\nu t=0$, therefore
$(x)\nu s+y$ is in the image of $w$: There is $z\in\tau_\Lambda C$
such that $(z)w=(x)\nu s+y$ and hence 
$$(z,x)\left[\smallmatrix 0 & w\cr 1 & -\nu s\endsmallmatrix\right]\;=\;
     (x, (z)w-(x)\nu s) \;=\; (x,y).$$
Comparing the two diagrams labelled $(*)$ and $(**)$ we obtain the 
following isomorphism.
$$\hobject {\tau_\Lambda B}X\gamma \quad\cong\quad 
  \hobject{\tau_\Lambda B}{\tau_\Lambda C\oplus \nu Q}{[\tau_\Lambda (g)
  \;v]} $$

\smallskip\noindent{\it Step 5}.
Write $\tau_{\Cal H}(g)=(h\: E\to F)$, so 
$\tau_{\Cal S}(f)=\Ker \tau_{\Cal H}(g)=(j\:D\to E)$ where 
$j=\Ker(h)$. 
Thus the object given by 
$\gamma = \Ker (\tau_{\Cal H}(g)\oplus(1_{\nu S}\:\nu S\to \nu S))$
in the previous step has the form
$$\hobject{D}E{j}\oplus\hobject 0{\nu S}{}
\quad\cong\quad \hobject{\tau_\Lambda B}X\gamma 
\quad\cong
  \hobject{\tau_\Lambda B}{\tau_\Lambda C\oplus \nu Q}{[\tau_\Lambda (g)
  \;v]},$$
where the second isomorphism has been established in Step 4.
By Claim 2 in the previous section, 
the right hand side is isomorphic to $\Mimo\tau_\Lambda(g)$,
up to an injective direct summand. Since neither $j$ nor 
$\Mimo\tau_\Lambda(g)$ has an injective direct summand 
--- recall that $j$
is a $\tau_{\Cal S}$-translate --- the Krull-Remak-Schmidt theorem in
$\Cal S(\Lambda)$ implies that 
$\tau_{\Cal S}(f)=j=\Mimo\tau_\Lambda(g)$. \qed

\bigskip
There is also the following dual version.

\medskip\noindent{\bf Theorem 5.2.} 
{\rm (The Auslander-Reiten translation in $\Cal F(\Lambda)$)}
\newline{\it 
For $\big(A\sto gB\big)$ an object in $\Cal F(\Lambda)$ the 
Auslander-Reiten translate is given as
$$\tau^-_{\Cal F}(g)\;=\;\Mepi\tau^-_\Lambda\Ker(g).$$
}\qed

\medskip\noindent
{\it Example.}  Let $\Lambda$ be a uniserial ring with maximal ideal
$\m$, Loewy length $n$ and socle $k=\m^{n-1}$, 
as in the Example in Chapter 1.  
Recall that the projective-injective indecomposable $(0\to\Lambda)$
in the category $\Cal S(\Lambda)$
has as sink map the inclusion $(0\to \m)\to(0\to\Lambda)$ and
as source map the inclusion $(0\to \Lambda)\to(k\to\Lambda)$, so
$\tau_{\Cal S}(k\to\Lambda)=(0\to\m)$.
We illustrate the formula in 5.1 by computing the powers of the 
Auslander-Reiten translation
for the module $(k\to\Lambda)$. First,
$$\tau_{\Cal S}\Big(\hobject k\Lambda{}\Big) 
        = \Mimo\tau_\Lambda\Cok\Big(\hobject k\Lambda{}\Big) 
        = \Mimo\tau_\Lambda\Big(\hobject\Lambda{\Lambda/k}{}\Big)
        = \Mimo\Big(\hobject0{\Lambda/k}{}\Big)
        = \Big(\hobject0{\Lambda/k}{}\Big)
        \cong \Big(\hobject0\m{}\Big)$$
confirming the above result.  Further translates are computed easily:
$$\eqalign{
\tau_{\Cal S}\Big(\hobject0\m{}\Big) 
        &=\Mimo\tau_\Lambda\Cok\Big(\hobject0\m{}\Big)
        = \Mimo\tau_\Lambda\Big(\hobject\m\m{}\Big)
        = \Big(\hobject\m\m{}\Big)                        \cr
\tau_{\Cal S}\Big(\hobject\m\m{}\Big) 
        &= \Mimo\tau_\Lambda\Big(\hobject\m0{}\Big)
        = \Mimo\Big(\hobject\m0{}\Big)
        = \Big(\hobject\m\Lambda{}\Big)                        \cr
\tau_{\Cal S}\Big(\hobject\m\Lambda{}\Big) 
        &= \Mimo\tau_\Lambda\Big(\hobject\Lambda{\Lambda/\m}{}\Big)
        = \Mimo\Big(\hobject0{\Lambda/\m}{}\Big)
        \cong \Big(\hobject0k{}\Big)                        \cr
\tau_{\Cal S}\Big(\hobject0k{}\Big) 
        &= \Mimo\tau_\Lambda\Big(\hobject kk{}\Big)
        = \Big(\hobject kk{}\Big)                        \cr
\tau_{\Cal S}\Big(\hobject kk{}\Big) 
        &= \Mimo\tau_\Lambda\Big(\hobject k0{}\Big)
        = \Big(\hobject k\Lambda{}\Big)                        \cr
}$$
So after six steps we are back at where we started.  This is not
a coincidence, as we will see in the next chapter.

\medskip\noindent{\it Definition:\/} Let 
$\Cal S(\Lambda)_{\Cal I}$ be the full subcategory of 
$\Cal S(\Lambda)$ consisting 
of all objects which have no nonzero injective direct summands.

\medskip\noindent{\bf Corollary 5.3.} {\it 
Every object in $\Cal S(\Lambda)_{\Cal I}$ has the form $\Mimo(f)$ for some 
morphism $f\:A\to B$ where the $\Lambda$-modules 
$A$ and $B$ have no nonzero injective direct summands.} \qed

\medskip\noindent{\it Definition:\/} By $\Cal H'(\Lambda)$ we denote
the morphism category for $\underlamod$.  Thus, the objects are 
the morphisms $(\bar f\:A\to B)$ in $\underlamod$, and a morphism 
in $\Cal H'(\Lambda)$ from $(\bar f\:A\to B)$ to $(\bar{f'}\:A'\to B')$
is given by a pair $(\bar a,\bar b)$ where $\bar a\:A\to A'$, $\bar b\:B'\to B$
are morphisms such that $\bar f\bar b=\bar a\bar {f'}$ holds in
$\underlamod$.

\medskip\noindent{\bf Corollary 5.4.} {\it 
The functor
$$F\:\Cal S(\Lambda)_{\Cal I}\longrightarrow \Cal H'(\Lambda),
        \qquad f\mapsto \bar f,$$
is dense, preserves indecomposables, and reflects isomorphisms.  
Thus, $F$ is a representation equivalence. 
}

\smallskip\noindent{\it Proof:\/} Clearly, $F$ is dense:
If the object $\bar f$ in $\Cal H'(\Lambda)$ is represented by $f$, then 
$F\Mimo(f)$ and $\bar f$ are isomorphic in $\Cal H'(\Lambda)$.

The functor $F$ preserves indecomposables:  Let $a$ be an 
indecomposable noninjective object in $\Cal S(\Lambda)$, then
$a$ occurs as the first term of an Auslander-Reiten sequence and hence has
the form $a=\tau_{\Cal S}(c)$ for some indecomposable nonprojective
object $c\in\Cal S(\Lambda)$. In particular, 
$F(a)=F\Mimo\big(\tau_\Lambda\Cok(c)\big)$ is isomorphic to
$\tau_\Lambda\Cok(c)$ in $\Cal H'(\Lambda)$ which is an indecomposable
object since $\Mimo$ is additive. 

To show that $F$ reflects isomorphisms, let two objects $f,g\in\Cal S(\Lambda)$
be given such that $F(f)$ and $F(g)$ are isomorphic in $\Cal H'(\Lambda)$.
By Corollary 5.3, $f=\Mimo(f_1)$ and $g=\Mimo(g_1)$ where
$f_1\:A_1\to B_1$ and $g_1\:C_1\to D_1$ 
are maps between $\Lambda$-modules with no nonzero injective direct summands. 
Since $\bar{f_1}$ and $\bar{g_1}$ are isomorphic in $\Cal H'(\Lambda)$
there are isomorphisms of $\Lambda$-modules $u\:A_1\to C_1$ and $v\:B_1\to D_1$
such that $f_1v-ug_1$ factorizes through an injective $\Lambda$-module.  
By Theorem 4.2, $f=\Mimo(f_1)$ and $g=\Mimo(g_1)$ are isomorphic. \qed

\medskip One can go a little bit further. 

\medskip\noindent{\it Definition:} Let $\Cal I$ be the
ideal in the category $\Cal S(\Lambda)$ of all morphisms which factorize 
through an injective object in $\Cal S(\Lambda)$. By 
$\Cal S(\Lambda)/\Cal I$ we denote the factor category of $\Cal S(\Lambda)$
modulo $\Cal I$.  Clearly, the functor
$F\:\Cal S(\Lambda)\to \Cal H'(\Lambda), f\mapsto \bar f$, annihilates
every morphism in $\Cal I$ and hence defines a functor 
$\bar F\:\Cal S(\Lambda)/\Cal I\to \Cal H'(\Lambda)$. 

\medskip\noindent{\bf Lemma 5.5.} {\it The functor $\bar F$ is full and dense.}

\smallskip\noindent{\it Proof:\/} Since the dense functor $F$ 
from Corollary~5.4 factorizes over $\bar F$, also $\bar F$ is dense.
To show that $\bar F$ is full, let $f\:A\to B$ and $g\:C\to D$ be 
objects in $\Cal S(\Lambda)$, let $e\:A\to I$ be an injective envelope
and $e'\:B\to I$ an extension of $e$. 
If $(\bar u,\bar v)$ is a morphism in $\Hom_{\Cal H'}(\bar f,\bar g)$
then there are maps $u\:A\to C$, $v\:B\to D$, $w\:I\to D$ such that 
$ug-fv=ew$. The following diagram is commutative
$$\CD A @>f>> B \cr @| @VV[1\;e']V \cr A @>[f\;e]>> B\oplus I \cr
        @VuVV  @VV\big[{v\atop w}\big]V \cr C @>>g> D \endCD $$
and hence represents two morphisms 
$\CD f@>\big(1_A,[1_B\;e']\big)>\phantom{)}> [f\;e] @>\big(u,\big[{u\atop w}\big]\big)>>g
\endCD$ in the category $\Cal S(\Lambda)$.  (Note that the 
first is a split monomorphism, so that $f$ and $[f\;e]$ become isomorphic
objects in the factor category $\Cal S(\Lambda)/\Cal I$.)
By applying the functor $F$ to the composition of the two morphisms, we
obtain $(\bar u,\bar v)$. \qed

\medskip\noindent{\it Remark:\/} The following example shows that the
functor $\bar F$ is not a categorical equivalence in general.
Let $\Lambda$ be a uniserial ring of Loewy length 3 and $m$ a generator of the
maximal ideal $\m$. Let $f$ be the object in $\Cal S(\Lambda)$ given by
the inclusion $f\:\m/\m^2\to \Lambda/\m^2$.
Then the multiplication $\mu_m$ by $m$ 
is a nonzero nilpotent endomorphism of $f$ which does not factorize
through an injective object in $\Cal S(\Lambda)$, but for which $F(\mu_m)=0$
holds. Thus, the functor $\bar F$ is not faithful.

        \bigskip
\centerline{\bf  6. Morphisms in the Stable Category}
       \nopagebreak \medskip
In this section we assume that $\Lambda$ is a self-injective algebra.  Then
the stable category $\underlamod$
modulo all morphisms which factorize through a 
projective-injective object, 
together with the suspension functor $\Omega^{-1}$,
becomes a triangulated category.  We recall the 
construction of triangles in this category, observe that the assignment
$f\mapsto \Cok\Mimo(f)$ for a morphism $f\:A\to B$ 
is related to the rotation of a triangle
in the stable category, and retrieve the formula 
$\tau_{\Cal S}^6(f)\cong f$
for the case that $\Lambda$ is a uniserial algebra.

\smallskip We recall from [H, Chapter I] that 
the stable category 
$\underlamod$ has the following standard triangles.
For a morphism $f\:A\to B$ in $\lamod$, take 
the short exact sequence given by the inclusion of $A$
in its injective envelope, $i\:A\to \I(A)$, and the
cokernel map $c\:\I(A)\to \Omega^{-1}A$, and form the
push-out diagram along $f$:
$$\CD   0 @>>> A @>i>> \I(A) @>c>> \Omega^{-1}A @>>> 0 \cr
        @. @VfVV  @VVjV @| \cr
        0 @>>> B @>>g> C @>>h> \Omega^{-1}A @>>> 0 \endCD$$
Then the standard triangle corresponding to $f$ is
$$ T(f)\: \qquad A \lto{\bar f} B\lto{\bar g} C \lto{\bar h} 
                \Omega^{-1}A$$
and each triangle in the stable category is isomorphic to a standard one
[H, Theorem I.2.6].

\smallskip A morphism between triangles $T\:A\to B\to C\to \Omega^{-1}A$ 
and $T'\: A'\to B'\to C'\to \Omega^{-1}A'$ 
consists of morphisms $\bar a$, $\bar b$, $\bar c$ in the stable category
which make the following diagram commutative.
$$\CD   A @>\bar f>> B @>>> C @>>> \Omega^{-1}A \cr
        @V\bar aVV @V\bar bVV @V\bar cVV @VV\Omega^{-1}\bar aV \cr
        A' @>>\bar{f'}> B' @>>\bar{g'}> C'@>>> \Omega^{-1}A'
        \endCD$$
The two triangles are isomorphic if $\bar a$,
$\bar b$, and $\bar c$ are isomorphisms.
From [H, I.1.2 and I.1.6] we recall that for the existence 
of a morphism between the triangles it is sufficient to have a map
$b\: B\to B'$ such that $\overline{fbg'}=0$.  Moreover, a morphism
$(\bar a, \bar b, \bar c)$ is an isomorphism if $\bar a$ and $\bar b$ 
are isomorphisms in $\underlamod$. Thus, isomorphisms between 
triangles are obtained from pairs $(\bar a,\bar b)$ 
of isomorphisms which make the left square in the above 
diagram commutative.
This is to say that if $(\bar a,\bar b)\: \big(A\sto{\bar f} B\big)\to
        \big(A'\sto{\bar{f'}}B'\big)$ is an isomorphism in the category
$\Cal H'(\Lambda)$ of morphisms in the stable category 
$\underlamod$, then the triangles $T(f)$ and $T(f')$ are isomorphic. 

\medskip Given a triangle
$T\: A\sto{\bar f} B\sto{\bar g}C\sto{\bar h}\Omega^{-1}A$, then the 
{\it rotation\/} of $T$ below also is a triangle, as required by 
the axioms of a triangulated category.
$$T^\R\:\qquad 
        B\sto{\bar g}C\sto{\bar h}\Omega^{-1}A\lto{-\Omega^{-1}\bar f}
                \Omega^{-1}B$$
This operation $T(f)\mapsto T(f)^\R$ in the triangulated 
category gives rise to a self\-equi\-valence $\bar f\mapsto {\bar f}^\R$ 
on the category $\Cal H'(\Lambda)$ of morphisms in the stable category.  
According to the following 
Lemma, a map $g$ in $\lamod$ representing  ${\bar f}^\R$ 
in $\underlamod$ is obtained as $g=\Cok\Mimo f$.

\medskip\noindent{\bf  Lemma 6.1.} {\it 
For a map $f\:A\to B$ in $\lamod$,
the two morphisms ${\bar f}^\R$ and $\overline{\Cok\Mimo f}$ are isomorphic 
in $\Cal H'(\Lambda)$.}

\smallskip\noindent{\it Proof:\/} In the following diagram with
exact rows, the top row
defines the cokernel map for $\Mimo(f)=[f\;e]$ 
while the third row is given by the
push-out diagram defining $T(f)$. 
$$\CD   0 @>>> A @>[f\;e]>> B\oplus \I\Ker (f)@>\Cok>>C' @>>> 0 \cr
        @. @| @V{\big[{1\atop0}{0\atop1}{0\atop0}\big]}VV @V[1\;0]VV \cr
        0 @>>> A @>[f\;e\;0]>> B\oplus \I\Ker (f)\oplus I' @>>> C'\oplus I' @>>> 0 \cr
        @. @|  @VdV\cong V @Vd'V\cong V \cr
        0@>>> A @>[f\;i]>> B\oplus \I(A)@>\big[{g\atop -j}\big]>> C @>>> 0 \cr
        @. @. @V{\big[{1\atop0}\big]}VV @| \cr
        @. @. B @>g>> C \endCD$$
All squares with the possible exception of the square at the
bottom are commutative:
Since the map $[f\;i]$ is a monomorphism 
it follows from Claim~2 in Chapter~4 that there is an injective module $I'$
and an isomorphism $d$ such
that the left square between the second and third row is commutative;
if $d'$ is the cokernel map then both squares are commutative.
Here the exact sequence in the first row is a direct summand of the 
sequence in the second row, the complement being the sequence
$0\to 0\to I'\sto1I'\to 0$. Finally, the map $g$ which represents the
second map in the triangle $T(f)$ occurs as a restriction of the
cokernel map in the third row.  

Note that all vertical maps become isomorphisms when considered 
in the stable 
category.  Thus, the two morphisms $\overline\Cok\:B\oplus\I\Ker (f)\to C'$
and $\bar g\:B\to C$ are isomorphic when considered as objects
in $\Cal H'(\Lambda)$.  
\qed

\medskip Since $\Lambda$ is a self-injective algebra, there is a third
selfequivalence on the morphism category $\Cal H'(\Lambda)$ 
(besides the suspension and the rotation) given by
the Auslander-Reiten translation $\tau_\Lambda$. 
According to [ARS, Proposition IV.3.7], the functors
$\tau_\Lambda$ and $\Cal N\Omega^2$ from $\underlamod$ to 
$\underlamod$ are isomorphic, where $\Cal N=D\Hom_\Lambda(-,\Lambda)$
is the Nakayama automorphism.
It follows that the 
Auslander-Reiten translation $\tau_\Lambda$ preserves triangles.
As a consequence, the functor $\tau_\Lambda$ commutes with the rotation 
$\bar f\mapsto {\bar f}^\R$, and also with the suspension $\Omega^{-1}$
in the sense that for each morphism $f\:A\to B$ 
there is a commutative diagram in which the vertical maps are isomorphisms.
$$\CD \tau_\Lambda\Omega^{-1}A @>\tau_\Lambda\Omega^{-1}f>>
                \tau_\Lambda\Omega^{-1}B \cr
        @V\eta'(f)V\cong V @V\eta''(f)V\cong V \cr
        \Omega^{-1}\tau_\Lambda A @>>\Omega^{-1}\tau_\Lambda f>
                \Omega^{-1}\tau_\Lambda B \endCD $$

\medskip
The selfequivalence on the triangulated category given by the 
rotation $T\mapsto T^\R$ yields isomorphisms
$T^{3\R}\cong -\Omega^{-1}T$ and $T^{6\R}\cong \Omega^{-2}T$.
We obtain the following consequence for the Auslander-Reiten translation
$\tau_{\Cal S}$ in the submodule category.

\medskip\noindent{\bf Theorem 6.2.} {\it
Suppose $\Lambda$ is a self-injective algebra.
If $(f\: A\to B)$ is an indecomposable nonprojective object in 
$\Cal S(\Lambda)$ then  there are the following isomorphisms in the 
morphism category $\Cal H'(\Lambda)$.}
$$\eqalign{(1)& \qquad\overline{\tau_{\Cal S}(f)} \;\cong\;
                        \tau_\Lambda\big(\overline{\Cok f}\big),\cr
  (2)& \qquad\overline{\tau_{\Cal S}^3(f)} \;\cong\;
                        -\tau_\Lambda^3\Omega^{-1}(\bar f), \;\text{and} \cr
  (3)& \qquad\overline{\tau_{\Cal S}^6(f)} \;\cong\;
                        \tau_\Lambda^6\Omega^{-2}(\bar f).} $$

\medskip\noindent{\bf Corollary 6.3.} {\it
Under the assumptions of the theorem, there are the following
isomorphisms in the submodule category $\Cal S(\Lambda)$:
$$\tau_S^3(f)\;\cong\;-\Mimo\tau_\Lambda^3\Omega^{-1}(f),\qquad
  \tau_S^6(f)\;\cong\;\Mimo\tau_\Lambda^6\Omega^{-2}(f)$$
}
\smallskip\noindent{\it Proof of Corollary 6.3:\/}
The functor $\Cal S(\Lambda)_I\to \Cal H'(\Lambda)$,
$f\mapsto \bar f$ in Corollary 5.3 reflects isomorphisms,
so the assertion follows from Theorem 6.2. \qed

\smallskip\noindent{\it Proof of Theorem 6.2.:\/} The first statement (1) follows from
Theorem 5.1 since any map $g$ is stably equivalent to $\Mimo (g)$.
For the proof of assertion (2) we use Theorem 5.1 to compute
$$
 \tau_{\Cal S}^3\big(A\sto fB\big)\;=\;\Mimo\tau_\Lambda\Cok
        \Mimo\tau_\Lambda\Cok\Mimo\tau_\Lambda\Cok(f).
$$
Then we obtain the following isomorphisms of objects in 
$\Cal H'(\Lambda)$:
$$\eqalign{ \overline{\tau_{\Cal S}^3(f)}\; 
        &\cong\; \tau_\Lambda\big(\overline{\Cok\Mimo\tau_\Lambda\Cok \Mimo
                \tau_\Lambda\Cok f}\big) \cr
        &\cong\; \tau_\Lambda\big(
                \overline{\tau_\Lambda\Cok\Mimo\tau_\Lambda\Cok f}^\R\big) \cr
        &\cong\; \tau_\Lambda^2 \big(
                \overline{\Cok\Mimo\tau_\Lambda\Cok f}^\R\big) \cr
        &\cong\; \tau_\Lambda^2 \big(
                \overline{\tau_\Lambda\Cok f}^{2\R} \big) \cr
        &\cong\; \tau_\Lambda^3 \big( \overline{\Cok f}^{2\R} \big) \cr
        &\cong\; \tau_\Lambda^3 \big( {\bar f}^{3\R} \big) \cr
        &\cong\; \tau_\Lambda^3 \big(-\Omega^{-1}(\bar f)\;\big)
}$$
where the first isomorphism is justified by (1),
the second, fourth, and sixth isomorphisms follow from Lemma 6.1, the third and
fifth equalities are come from
the commutativity of $\tau_\Lambda$ with the rotation,
and the last map is an isomorphism since a threefold rotation of a triangle is 
obtained by applying the functor $-\Omega^{-1}$. 

In order to deduce the third assertion from the second, pick a representative
map $g\:\Omega^{-1}A\to\Omega^{-1}B$ 
in $\lamod$ for the morphism $\overline{\tau_{\Cal S}^3(f)}$
in $\underlamod$ such that $\Omega^{-1}A$ and $\Omega^{-1}B$ have no
nonzero injective direct summands.  
Then $\Mimo(g)$ is an indecomposable nonprojective object in $\Cal S(\Lambda)$
and by Proposition 4.1 its isomorphism 
class does not depend on the choice of the map $g$. 
The following morphisms in $\underlamod$ 
are isomorphic objects in $\Cal H'(\Lambda)$.
$$\overline{\tau_{\Cal S}^6(f)} 
        \;\cong\; \overline{\tau_{\Cal S}^3\Mimo (g)}
        \;\cong\; -\tau_\Lambda^3\Omega^{-1}(\bar g)
        \;\cong\; -\tau_\Lambda^3\Omega^{-1}
                \big(-\tau_\Lambda^3\Omega^{-1}(\bar f)\big) 
        \;\cong\; \tau_\Lambda^6\Omega^{-2}(\bar f).
$$
\qed

\medskip

We conclude this section with three applications.  

\medskip\noindent{\bf Corollary 6.4.}  {\it
Suppose $\Lambda$ is a self-injective algebra such that 
$\tau_\Lambda$ coincides with $\Omega^2$. 
If $(f\: A\to B)$ is an indecomposable nonprojective object in 
$\Cal S(\Lambda)$, then there is an isomorphism of objects
$$
     \tau_{\Cal S}^3(f) \cong
                        -\Mimo\Omega^{5} (f)
$$
in $\Cal S(\Lambda)$.}

\smallskip\noindent{\it Proof:\/} Since $\tau_\Lambda$ 
coincides with $\Omega^2$, 
we can simplify the expression in formula 2 of Theorem 6.2 and see that 
$$
    \overline{\tau_{\Cal S}^3(f)} \cong
                        -\tau_\Lambda^3\Omega^{-1}(\bar f)   \cong
                        -\Omega^6\Omega^{-1}(\bar f)  \cong
                        -\Omega^5(\bar f)
$$ 
in $\Cal H'(\Lambda)$. The functor 
$\Cal S(\Lambda)_{\Cal I}\to \Cal H'(\Lambda)$, $f\mapsto \bar f$,
in Corollary 5.3 reflects isomorphisms, so 
$\tau_{\Cal S}^3(f)$ and $-\Omega^5(f)$
are isomorphic in $\Cal S(\Lambda)$. \qed

\medskip
Note that for any symmetric algebra, the functors  	
$\tau_\Lambda$ and $\Omega^2$ coincide (see for example [ARS, Proposition IV.3.8]), thus we
can apply Corollary 6.4 in this case.

\medskip\noindent{\bf Corollary 6.5.}  {\it
Let $\Lambda$ be a commutative uniserial algebra.
Then for an indecomposable nonprojective 
object $(f\:A \to B)$ in $\Cal S(\Lambda)$, there is an isomorphism of objects
$$
 \tau_{\Cal S}^6(f)\;\cong\; f
$$
in $\Cal S(\Lambda).$}

\smallskip\noindent{\it Proof:\/}
Since $\Lambda$
is a commutative uniserial algebra, 
all the  functors $\tau_\Lambda,$ $\Omega^{2}$ and $\Omega^{-2}$ are 
equivalent to the identity functor on $\underlamod$,
thus Corollary 6.4 shows that
$\tau_{\Cal S}^6(f)$ and $\Omega^{10}(f)\; \cong\; f$ are isomorphic objects
of $\Cal S(\Lambda)$. \qed

\medskip\noindent{\it Definition:\/} By $\Bbb A_\infty^\infty$ we denote
the doubly infinite linear quiver 
$$\cdots \quad \lfrom\alpha \bullet^{-1} \lfrom\alpha\bullet^0
        \lfrom\alpha\bullet^1\lfrom\alpha \quad \cdots $$
The path algebra $k\Bbb A_\infty^\infty$ 
of this quiver is the associative $k$-algebra with basis the paths
in $\Bbb A_\infty^\infty$. If $\alpha^n$ denotes the ideal spanned by
all paths of length at least $n$, then 
the factor algebra $\Lambda=k\Bbb A_\infty^\infty/\alpha^n$
is a locally bounded associative $k$-algebra.  
A $\Lambda$-module $A$ 
consists of a sequence $(A_i)_{i\in\Bbb Z}$ of $k$-modules 
together with a sequence $(\alpha_i\:A_i\to A_{i-1})_{i\in\Bbb Z}$
of linear maps.
By $A[\ell]$ we denote the {\it shifted\/} module given by the spaces
$(A_{i-\ell})_i$ and the maps $(\alpha_{i-\ell})_i$.

\medskip\noindent{\bf Corollary  6.6.}  {\it 
Let $\Lambda$ be the 
associative algebra $k\Bbb A_\infty^\infty/\alpha^n$ where $k$ is a field.
For an indecomposable nonprojective 
object $(f\:A\to B)$ in $\Cal S(\Lambda)$, the following formula holds.}
$$\tau_{\Cal S}^6(f)\;\cong\; f[n-6]$$

\smallskip\noindent{\it Proof:\/}
The Auslander-Reiten translation $\tau_\Lambda$ is given by the shift
$A\mapsto A[-1]$ along the arrow $\alpha$, hence
the functor $\tau_\Lambda$ on the stable category
$\underlamod$ preserves triangles.
Also, $\tau_\Lambda$ commutes with
$\Omega^{-1}$ and with the rotation in a triangle. 
Moreover, for a nonprojective
indecomposable $\Lambda$-module $A$, the process of taking the cokernel of
the injective envelope twice yields the module $A[n]$.
With these adjustments, 
the claim follows from (3) in Theorem~6.2 
as in the proof of Corollary~6.4. \qed

        \bigskip
\centerline{\bf  7. Auslander-Reiten Sequences}
\nopagebreak        \medskip
In this section we show that ``most'' Auslander-Reiten sequences in the
category $\Cal S(\Lambda)$ become split exact sequences in
the category $\lamod$, when restricted to the 
short exact sequence of the submodules, or to the short exact sequence of
the big modules. We describe the exceptions in detail.
The remaining sink and source maps are associated with the
projective and the injective objects and have been specified
in Chapter~1.\ We only need to assume here that 
$\Lambda$ is a locally bounded associative algebra. 
First we deal with the exceptions.

\medskip\noindent
{\bf Proposition 7.1.} {\rm (Auslander-Reiten sequences with
        components not split exact)} \newline
{\it 
Let $0\to A\sto fB\sto gC\to 0$ be an Auslander-Reiten sequence in $\lamod$.
\smallskip
\item{1.} The Auslander-Reiten sequence in $\Cal S(\Lambda)$ ending at
$(0\to C)$ has the form
$$\CD 0 @>>> \hobject AA1 @>1>f> \hobject ABf @>0>g> \hobject 0C{} @>>> 0.
\endCD $$
\item{2.} With $e\:A\to \I(A)$ an injective envelope, the Auslander-Reiten
sequence in $\Cal S(\Lambda)$ ending at $(1_C\: C\to C)$ has the form
$$\CD 0 @>>> \hobject A{\I(A)}e @>f>[1\;0]> \hobject B{\I(A)\oplus C}b
        @>g>{\left[{0\atop 1}\right]}> \hobject CC1 @>>> 0
\endCD $$
where $b$ is the map $[e'\;g]$ with 
$e'\:B\to\I(A)$ an extension of the map $e$.

}

\smallskip\noindent
{\it Proof:\/}
1. We show that the map $(0,g)$ is minimal right almost split.
Clearly, this map is right minimal and not a split epimorphism.
In order to show that $(0,g)\: \big(A\sto fB\big) \to \big(0\to C\big)$
is right almost split, 
let $(t',t):(x'\:X'\to X)\to(0\to C)$ be a test map which is not a split
epimorphism.  Then $t\: X\to C$ is not a split epimorphism in $\lamod$,
so there is $u\:X\to B$ such that $t=ug$. 
Since the composition $X'\sto{x'}X\sto uB\sto gC$ is zero,
there is $u': X'\to A$ such that $x'u=u'f$.
Thus, $(u',u)$ is a morphism which satisfies $(t',t)=(u',u)(0,g)$.

\smallskip
2. A straightforward argument shows that the map $(g,\left[{0\atop 1}\right])$
is right minimal. 
We verify that a test map $(t',t)\: (x'\:X'\to X)\to(1_C\:C\to C)$
which is not a split epimorphism factorizes over $(g,\left[{0\atop 1}\right])$.
Since $t'\:X'\to C$ is not a split epimorphism, there is
$u'\: X'\to B$ such that $t'=u'g$. 
Let $u_1\:X\to \I(A)$ be an extension of $u'e'\:X'\to \I(A)$
to $X$ and put $u=[u_1\;t]\:X\to \I(A)\oplus C$.
Then $(u',u)$ is  a morphism in $\Cal S(\Lambda)$ such that
$(t',t)=(u',u)(g,\left[{0\atop 1}\right])$. \qed

\medskip
The remaining Auslander-Reiten sequences in $\Cal S(\Lambda)$
are made up from two split exact sequences.

\medskip\noindent
{\bf Proposition 7.2.} {\rm (Auslander-Reiten sequences with components split
exact)} \newline {\it
Suppose that $(c:C'\to C)$ is an indecomposable object in $\Cal S(\Lambda)$
such that the  morphism $c$ in $\lamod$ is not split monic. If
$$\CD 0 @>>> \hobject{A'}A{a} @>f'>f> \hobject{B'}B{b} @>g'>g>
             \hobject{C'}C{c} @>>> 0 
\endCD $$
is an Auslander-Reiten sequence in $\Cal S(\Lambda)$ then both sequences
in $\lamod$,
$$0\to A'\sto{f'} B' \sto{g'} C' \to 0\T{and}0\to A\sto fB\sto gC\to 0$$
are split exact.
}

\smallskip\noindent
{\it Proof:\/} Since the map $c$ is not a split monomorphism, 
the test maps
$$(0,1)\:(0\to C)\to (C'\sto{c} C)\quad\text{and}\quad
  (1,c)\:(C'\sto 1 C') \to (C'\sto{c}C)$$
are not split epimorphisms and hence factorize over $(g',g)$. 
Thus, both $g$ and $g'$ are split epimorphisms. \qed

\medskip We combine this result with Theorem 5.1.

\medskip\noindent
{\bf Corollary 7.3.} {\rm (The middle term of an Auslander-Reiten sequence)}
\newline {\it
Let $(c:C'\to C)$ be an indecomposable object in $\Cal S(\Lambda)$
such that the  morphism $c$ in $\lamod$ is not split monic. 

\smallskip\item{1.} 
The Auslander-Reiten sequence ending in $(c\:C'\to C)$ has the form
$$\CD 0 @>>> \hobject{A'}Aa
                @>\big[{1\atop 0}\big]>\big[{1\atop0}\big]> 
        \hobject{A'\oplus C'}{A\oplus C}{b} 
                @>[0\;1]>[0\;1]>
             \hobject{C'}C{c} @>>> 0 
\endCD $$

\item{2.} The map $b$ defining the middle term is given as follows by
a map $h\: C'\to A$.
$$b=\pmatrix a & 0 \cr h & c \endpmatrix\: A'\oplus C'\to A\oplus C$$

\item{3.} The first term 
$\big(A'\sto aA\big)=\tau_{\Cal S}\big(C'\sto cC\big)$
of the Auslander-Reiten sequence ending in $(c\:C'\to C)$ is isomorphic
to $\Mimo\tau_\Lambda\Cok(c)$, in particular, $A'=\tau_\Lambda C$
and $A=\tau_\Lambda C''\oplus I$ where $I$ is an injective $\Lambda$-module
and $C''=\Cok(c)$.\qed

}

\medskip
To conclude this chapter, we state the 
following dual results for the category $\Cal F(\Lambda)$.

\medskip\noindent
{\bf Proposition 7.4.} {\rm (AR-sequences in $\Cal F(\Lambda)$ with
components not split exact)} \newline {\it 
Suppose that $0\to A\to B\to C\to 0$ is an Auslander-Reiten sequence
in $\lamod$.
\smallskip
\item{1.} The Auslander-Reiten sequence in $\Cal F(\Lambda)$ 
starting at $(A\to 0)$ has the form
$$\CD 0 @>>> \hobject A0{} @>f>0> \hobject BCg @>g>1> \hobject CC1 @>>> 0.
\endCD $$
\item{2.} The Auslander-Reiten sequence in $\Cal F(\Lambda)$
starting at $(1_A\:A\to A)$ has the form
$$\CD 0 @>>> \hobject AA1 @>{[1\;0]}>f> \hobject{A\oplus P}Bb
        @>{\left[{0\atop 1}\right]}>g> \hobject PC\pi @>>> 0
\endCD $$
where $\pi\: P\to C$ is a projective cover, and where the map $b$ has
the form $b=[f\;\hat\pi]^t$ with $\hat\pi\:P\to B$ a lifting of $\pi$
to $B$. \qed

}

\medskip\noindent
{\bf Proposition 7.5.} {\rm (AR-sequences in $\Cal F(\Lambda)$ with split
exact components)} \newline {\it
Suppose that $(a\:A\to A'')$ is indecomposable in $\Cal F(\Lambda)$
such that the epimorphism $a$ is not a split epimorphism. If
$$\CD 0@>>> \hobject A{A''}{} @>f>f''> \hobject B{B''}{} @>g>g''> 
            \hobject C{C''}{} @>>> 0
\endCD $$
is an Auslander-Reiten sequence in $\Cal F(\Lambda)$ then both
sequences in $\lamod$, 
$$0\to A\sto fB \sto gC \to 0\T{and}
0\to A''\sto{f''}B''\sto{g''}C''\to 0$$ 
are split exact.  \qed
}

\medskip\noindent
{\bf Corollary 7.6.} {\rm (The middle term of an AR-sequence in 
$\Cal F(\Lambda)$)} \newline {\it
Let $(a:A\to A'')$ be an indecomposable object in $\Cal F(\Lambda)$
such that the map $a$ in $\lamod$ is not a split epimorphism. 

\smallskip\item{1.} The Auslander-Reiten sequence starting at $(a\:A\to A'')$ 
has the form
$$\CD 0 @>>> \hobject{A}{A''}a
                @>\big[{1\atop 0}\big]>\big[{1\atop0}\big]> 
        \hobject{A\oplus C}{A''\oplus C''}{b} 
                @>[0\;1]>[0\;1]>
             \hobject{C}{C''}{c} @>>> 0 
\endCD $$

\item{2.} The morphism $b$ defining the middle term is given as follows 
by some map $h\: C\to A''$.
$$b=\pmatrix a & 0 \cr h & c \endpmatrix\: A\oplus C\to A''\oplus C''$$

\item{3.} The last term 
$\big(C\sto cC''\big)=\tau^-_{\Cal F}\big(A\sto aA''\big)$
of the Auslander-Reiten sequence starting in $(a\:A\to A'')$ is isomorphic
to $\Mepi\tau^-_\Lambda\Ker(a)$, in particular, $C=\tau^-_\Lambda A'\oplus P$
and $C''=\tau^-_\Lambda A$ where $P$ is a projective $\Lambda$-module
and $A'=\Ker(a)$. \qed

}

%
%
%
%
%

        \bigskip\bigskip
\frenchspacing

\noindent
{\bf References}     \nopagebreak   {\baselineskip=9pt \rmk
\parindent=1.5truecm
                                        \medskip\smallskip\noindent
 \item{[ARS]} M.\ Auslander, I\. Reiten, S.~O.~Smal\o: 
  {\itk Representation Theory of Artin Algebras.} 
        Cambridge University Press. (1995).

 \item{[AS]} M.\ Auslander, S.~O.~Smal\o: 
  {\itk Almost split sequences in subcategories,}
  Journal of Algebra,~{\bfk 69} (1981), 426-454.
        \smallskip
\item{[GR]} P.~Gabriel and A.~V.~Ro\u\i ter, 
        {\itk Representations of finite 
dimensional algebras (with a chapter by Bernhard Keller),}
in: Encyclopedia of mathematical Sciences {\bfk 73}, Algebra VIII,
1--177, Springer, Berlin 1992.

\smallskip\item{[H]} D.\ Happel, {\itk Triangulated Categories in the 
Representation Theory of Finite Dimensional Algebras,} London Mathematical
Society Lecture Notes series {\bfk 119}, ix+208pp, Cambridge University 
Press, Cambridge 1988.

\item{[RS1]} C.M\. Ringel, M\. Schmidmeier: 
{\itk Invariant subspaces of nilpotent linear operators,}
  In preparation.

\item{[RS2]} C.M\. Ringel, M\. Schmidmeier: 
{\itk Submodule categories of wild representation type,}
Journal for Pure and Applied Algebra (to appear), 1-12.

\item{[S]} A\. Skowro\'nski: 
{\itk Tame triangular matrix algebras over Nakayama algebras,}
J.~London Math.\ Soc.\ {\bfk 34} (1986), 245-264.
  
}
        \bigskip\medskip\noindent
{\rmk Claus Michael Ringel,
Fakult\"at f\"ur Mathematik, Universit\"at Bielefeld,
\par\noindent  POBox 100\,131, \ D-33\,501 Bielefeld 
\par\noindent {\ttk ringel\@math.uni-bielefeld.de}
        \medskip\noindent
Markus Schmidmeier, 
Department of Mathematical Sciences, Florida Atlantic University,
\par\noindent Boca Raton, Florida 33431-0991
\par\noindent {\ttk markus\@math.fau.edu}
}

\bye